\documentclass[11pt]{article}
\RequirePackage[OT1]{fontenc}
\RequirePackage[aos,amsthm,amsmath,natbib]{imsart}
\usepackage{amsfonts,graphicx,color}
\newcommand {\apgt} {\ {\raise-.5ex\hbox{$\buildrel>\over\sim$}}\ }
\newcommand {\aplt} {\ {\raise-.5ex\hbox{$\buildrel<\over\sim$}}\ }

\newcommand{\yv}{{\bf y}}

\newcommand{\epsv}{{\bf \epsilon}}
\newcommand{\zv}{{\bf z}}
\newcommand{\wv}{{\bf w}}

\begin{document}


\begin{frontmatter}

\title{
How Many Iterations are Sufficient for Semiparametric Estimation?}
\runtitle{$K$-step Semiparametric Estimation} \today
\author{Guang Cheng\thanksref{t1}}


\runauthor{Guang Cheng}

\affiliation{Purdue University}

\address{Guang Cheng\\
Department of Statistics\\
Purdue University\\
West Lafayette, IN 47907-2066\\Email: chengg@purdue.edu }

\thankstext{t1}{Guang Cheng is Assistant Professor,
Department of Statistics, Purdue University}

\begin{abstract}
A common practice in obtaining a semiparametric efficient estimate is through iteratively maximizing the (penalized) log-likelihood w.r.t. its Euclidean parameter and functional nuisance parameter via Newton-Raphson algorithm. The purpose of this paper is to provide a formula in calculating the minimal number of iterations $k^\ast$ needed to produce an efficient estimate $\widehat\theta_n^{(k^\ast)}$ from a theoretical point of view. We discover that (a) $k^\ast$ depends on the convergence rates of the initial estimate and nuisance estimate; (b) more than $k^\ast$ iterations, i.e., $k$, will only improve the higher order asymptotic efficiency of $\widehat\theta_n^{(k)}$; (c) $k^\ast$ iterations are also sufficient for recovering the estimation sparsity in high dimensional data. These general conclusions hold, in particular, when the nuisance parameter is not estimable
at root-n rate, and apply to semiparametric models estimated under various regularizations, e.g., kernel or penalized estimation. This paper provides a first general theoretical justification for the ``one-/two-step iteration" phenomena observed in the literature, and may be useful in reducing the bootstrap computational cost for the semiparametric models.

\end{abstract}


\end{frontmatter}

\newtheorem{theorem}{\indent \sc Theorem}
\newtheorem{corollary}{\indent \sc Corollary}
\newtheorem{lemma}{\indent \sc Lemma}
\newtheorem{proposition}{\indent \sc Proposition}
\newtheorem{remark}{\indent \sc Remark}
\newcommand{\phif}{\textsc{igf}}
\newcommand{\sign}{\mbox{sign}}
\newcommand{\phgf}{\textsc{gf}}
\newcommand{\fix}{$\textsc{gf}_0$}
\newcommand{\mb}[1]{\mbox{\bf #1}}
\newcommand{\Exp}[1]{\mbox{E}\left[#1\right]}
\newcommand{\pr}[1]{\mbox{P}\left[#1\right]}
\newcommand{\pp}[0]{\mathbb{P}}
\newcommand{\ee}[0]{\mbox{E}}
\newcommand{\re}[0]{\mathbb{R}}
\newcommand{\argmax}[0]{\mbox{argmax}}
\newcommand{\argmin}[0]{\mbox{argmin}}
\newcommand{\ind}[0]{\mbox{\Large\bf 1}}
\newcommand{\narrow}{\stackrel{n\rightarrow\infty}{\longrightarrow}
}
\newcommand{\weakpn}{\stackrel{P_n}{\leadsto}}
\newcommand{\weakpnboot}{\mbox{\raisebox{-1.5ex}{$\stackrel
{\mbox{\scriptsize $P_n$}}{\stackrel{\mbox{\normalsize
$\leadsto$}} {\mbox{\normalsize $\circ$}}}$}}\,}
\newcommand{\ol}[1]{\overline{#1}}
\newcommand{\avgse}[1] { \bar{\widehat{\Sigma}}_{#1} }
\newcommand{\mcse}[1]  { \Sigma^{*}_{#1} }
\newcommand{\po}{\textsc{po}}
\newcommand{\ph}{\textsc{ph}}
\newcommand{\x}{\mathbf{x}}
\newcommand{\xw}{\mathbf{xw}}

\def\boxit#1{\vbox{\hrule\hbox{\vrule\kern6pt
          \vbox{\kern6pt#1\kern6pt}\kern6pt\vrule}\hrule}}
\def\jhcomment#1{\vskip 2mm\boxit{\vskip 2mm{\color{red}\bf#1}
 {\color{blue}\bf -- JH\vskip 2mm}}\vskip 2mm}

\def\jhcommentin#1{{(\color{red}\bf{#1}{\color{blue}\bf -- JH})}}



\def\boxit#1{\vbox{\hrule\hbox{\vrule\kern6pt
          \vbox{\kern6pt#1\kern6pt}\kern6pt\vrule}\hrule}}

\def\chengcomment#1{\vskip 2mm\boxit{\vskip 2mm{\color{blue}\bf#1}
{\color{red}\bf -- Cheng\vskip 2mm}}\vskip 2mm}

\def\chengcommentin#1{{(\color{blue}\bf{#1}{\color{red}\bf -- Cheng}})}

\section{Introduction}
Semiparametric models indexed by a Euclidean parameter of interest $\theta\in\Theta\subset\mathbb{R}^d$ and an infinite-dimensional nuisance parameter $\eta\in\mathcal{H}$ are proven to be useful in a variety of contexts, e.g., \citep{h96,mrv97,fhw95,cfgw97,ss94,rcl96,ac03,lsv08, bmm09}. The semiparametric MLE for $\theta$ can be viewed as a solution of the implicitly defined efficient score function whose nonparametric estimation is only possible in some special cases, e.g., \cite{h96}. Therefore, it is generally hard to solve the MLE from the efficient score function analytically or numerically. A common practice is to maximize the log-profile likelihood
\begin{eqnarray}
\log pl_n(\theta)=\sup_{\eta\in\mathcal{H}}\log lik_n(\theta,\eta),\label{prolik}
\end{eqnarray}
where $lik_n(\theta,\eta)$ is the likelihood given $n$ data, via some optimization algorithm. For example, the Newton-Raphson algorithm is applied to the partial likelihood of the Cox model in the software $\mathbf{R}$ (with the command {\it coxph}).

A general algorithm of obtaining a semiparametric efficient estimate of $\theta$ is to iteratively maximize the log-likelihood w.r.t. $\theta$ and $\eta$ as follows:

{\it General Semiparametric Iterative Estimation Algorithm}
\begin{itemize}
\item[I.] Identify an initial estimate $\widehat\theta^{(0)}_n$;

\item[II.] Construct the corresponding nuisance estimate $\widehat\eta(\widehat\theta^{(0)}_n)$ either by pure nonparametric approach, e.g., isotonic estimation, or under some regularization, e.g., kernel or sieve estimation;

\item[III.] Apply the Newton-Raphson (NR) or other optimization algorithm to
    \begin{eqnarray}
    \widehat S_n(\theta)=\log lik_n(\theta,\widehat\eta(\theta)),\label{gpl-def}
    \end{eqnarray}
    at $\theta=\widehat\theta^{(0)}_n$ to obtain $\widehat\theta^{(1)}_n$;

\item[IV.] Repeat steps II-III $k^\ast$ iterations until $$|\widehat S_n(\widehat\theta^{(k^\ast)}_n)-\widehat S_n(\widehat\theta^{(k^\ast-1)}_n)|\leq\epsilon$$ for some pre-determined sufficiently small $\epsilon$.
\end{itemize}
Note that $\widehat S_n(\theta)$ defined in (\ref{gpl-def}) is also called the generalized profile likelihood in \cite{sw92}. If $\widehat\eta(\theta)$ is the nonparametric MLE (NPMLE) for any fixed $\theta$, then $\widehat S_n(\theta)$ is just the profile likelihood defined in (\ref{prolik}). The above likelihood estimation procedure or its M-estimation analog has been extensively implemented in the literature. Here is an incomplete list: (i) Odds-Rate Regression Model under Survival Data, e.g., \cite{h96,mrv97}; (ii) Semiparametric Regression under Shape Constraints, e.g., \cite{c09, bmm09}; (iii) Logistic Regression with Missing Covariates, e.g., \cite{rcl96}; (iv) Generalized Partly Linear (Single Index) Model, e.g., \cite{fhw95,cfgw97}; (v) Conditionally Parametric Model, e.g. \cite{sw92, ss94}; (vi) Semiparametric Transformation Model, e.g., \cite{lsv08}. In addition, the above iterative procedure can also be adapted to the penalized estimation and selection of the semiparametric models by using a different criterion function than (\ref{gpl-def}), see \cite{cz10,bm05,mvg97}. We will discuss that scenario in Section~\ref{penest}. However, in all the above papers, $k^\ast$ or $\epsilon$ is {\it arbitrarily} chosen in practice.

The main purpose of our paper is to answer ``How Many Iterations Do We Really Need in Semiparametric Estimation?" from a {\it theoretical} point of view. We provide a general formula in calculating the minimal number of iterations $k^\ast$ needed to produce a semiparametric efficient $\widehat\theta_n^{(k^\ast)}$. Specifically, we discover that (a) $k^\ast$ depends on the convergence rates of $\widehat\theta_n^{(0)}$ and $\widehat\eta(\theta)$; (b) more than $k^\ast$ iterations, i.e., $k$, will not change the limiting distribution of $\widehat\theta_n^{(k)}$, but will improve its higher order asymptotic efficiency; (c) $k^\ast$ iterations are also sufficient for recovering the estimation sparsity under high dimensional data. These general conclusions hold, in particular, when the nuisance parameter is not estimable at root-n rate, and apply to semiparametric models estimated under various regularizations, e.g., kernel or penalized estimation. Note that the convergence rate of the regularized estimate $\widehat\eta(\theta)$ is determined by the related smoothing parameters, e.g., the bandwidth order in kernel estimation. Moreover, our construction of the efficient estimate does not require knowing the form of the implicitly defined efficient score function or apply the sample splitting technique and the drop-one-out trick required in the classical literature, i.e., \citep{b82,s86,k87,s87}. A general strategy of identifying $\widehat\theta_n^{(0)}$ with proper convergence rate is also considered. The technical challenge of this paper is that $\widehat S_n(\theta)$ in practice may not have an explicit form or is not continuous/smooth.

As far as we are aware, our paper provides a first general theoretical justification for the ``one-/two-step iteration" phenomenon, i.e., $k^\ast=1, 2$, observed in the semiparametric literature. However, we find that more iterations are absolutely necessary if $\eta$ is estimated at a very slow rate. For example, we need 8 iterations to achieve the efficiency in conditionally exponential models, see Table 3. Moreover, our results are readily extended to the bootstrap estimation by combining with the most recent bootstrap consistency results obtained for semiparametric models in \cite{cz10}. Therefore, we expect to significantly reduce bootstrap computational cost, which is very high in semiparametric models, after knowing $k^\ast$ for each bootstrap sample. See \citep{a02} for similar ideas but applied to the parametric models. Due to the space limitation, we only consider the NR algorithm based on original sample in this paper, but notice that the extensions to the slight modifications of NR are possible by considering the discussions in Page 534 of \citep{r88}.

Section~\ref{secrev} provides some necessary background material on the semiparametric estimation. In Section~\ref{semisec}, we consider the semiparametric maximum likelihood estimation in which $\widehat S_n(\theta)$ is the possibly nonsmooth profile likelihood (\ref{prolik}). In Section~\ref{smre}, we consider the semiparametric estimation under two types of regularization, i.e., kernel estimation and penalized estimation, in which $\widehat S_n(\theta)$ is smooth. In that section, we also consider the sparse and efficient estimation of the partial linear models as an important application of penalized estimation. In Section~\ref{iniest}, we propose two grid search algorithms for identifying the initial estimate whose convergence rate will be rigorously proven. Several semiparametric models ranging from survival models, mixture models to conditionally exponential models are treated to illustrate the applicability of our theories. All the proofs are postponed to the Appendix.

\section{Preliminary}\label{secrev}
We assume that the data $X_1,\ldots,X_n$ are i.i.d. throughout the paper. In what follows, we first briefly review the concepts of the efficient score function and the least favorable curve (LFC), and then relate the estimation of LFC to that of $\theta$ as discussed in \cite{sw92}. Unless otherwise specified, the notation $E$ is reserved for the expectation taken under $(\theta_0,\eta_0)$.

The score functions for $\theta$ and $\eta$ are defined as,
respectively,
\begin{eqnarray}
\dot{\ell}_{0}(X_i)&=&\frac{\partial}{\partial\theta}\log
lik(X_i;\theta_0,\eta_0),\nonumber\\
A_{\theta_0,\eta_0}h(X_i)&=&\frac{\partial}{\partial t}|_{t=0}\log
lik(X_i;\theta_0,\eta(t)),\label{scoeta}
\end{eqnarray}
where $h$ is a ``direction" along which $\eta(t)\in\mathcal{H}$
approaches $\eta_0$  as $t\rightarrow 0$. $A_{\theta_0,\eta_0}: \mathbf{H}\mapsto L_2^0(P_{\theta_0,\eta_0})$ is the score operator for $\eta$, where $\mathbf{H}$ is some closed and linear diection set. The efficient score function $\widetilde{\ell}_0$ is
defined as the residual of the projection of
$\dot{\ell}_{0}$ onto the tangent space
$\mathcal{T}$, which is defined as the closed linear span of the tangent set $\{A_{\theta_0,\eta_0}H=(A_{\theta_0,\eta_0}h_1,\ldots,A_{\theta_0,\eta_0}
h_d)': h_j\in\mathbf{H}\}$. Therefore, we can write the efficient score function at $(\theta_0,\eta_0)$ as
\begin{eqnarray}
\widetilde{\ell}_0=\dot{\ell}_{0}-\Pi_{0}\dot\ell_{0},\label{effscor}
\end{eqnarray}
where $\Pi_{0}\dot\ell_{0}=\arg\min_{k \in\mathcal{T}}
E\|\dot{\ell}_{0}-k\|^2$. The variance of $\widetilde{\ell}_0$ is defined as the efficient
information matrix $\widetilde{I}_0$. The inverse of $\widetilde I_0$ is shown to be
Cram\'{e}r-Rao bound for estimating $\theta$ in the presence of an infinite dimensional $\eta$, see \cite{bkrw98}.

A main idea of estimating $\theta$ is to reduce a high dimensional
semiparametric model to a low dimensional random submodel of the
same dimension as $\theta$ called the least favorable
submodel (LFS). The LFS can be constructed as $t\mapsto \log
lik(t,\eta_\ast(t))$ and satisfies
\begin{eqnarray}
\eta_\ast(\theta_0)=\eta_0.\label{etatheta0}
\end{eqnarray}
and
\begin{eqnarray}
\frac{\partial}{\partial t}|_{t=\theta_0}\log
lik(t,\eta_\ast(t))=\widetilde{\ell}_0\label{lfs}
\end{eqnarray}
Note that the LFS may not exist unless $\Pi_0\dot\ell_0$ can be expressed as a nuisance score (the tangent set is closed). In all our examples, the LFS exists or can be approximated sufficiently closely. The $\eta_\ast(t)$ in the LFS is called as the least favorable curve. Under regularity conditions, it is shown that
\begin{eqnarray}
\eta_\ast(t)=\arg\sup_{\eta\in\mathcal{H}} E\log lik
(t,\eta)\;\;\;\mbox{for any fixed}\;t\in\Theta.\label{lfcurve}
\end{eqnarray}
By (\ref{lfcurve}) and standard arguments, we can establish that the maximizer of
$$S_n(\theta)\equiv \sum_{i=1}^{n}\log lik(\theta,\eta_\ast(\theta))(X_i)$$
is semiparametric efficient. In addition, based on (\ref{lfs}), we can derive that
\begin{eqnarray}
\;\;\;\;\;\widetilde{I}_0=E\left(\frac{\partial\log
lik(t,\eta_\ast(t))}{\partial t}|_{t=\theta_0}\right)^{\otimes 2}=
-E\left(\frac{\partial^2\log
lik(t,\eta_\ast(t))}{\partial t^2}|_{t= \theta_0}\right).\label{gpl-info2}
\end{eqnarray}

Recall that $\widehat S_n(\theta)=\sum_{i=1}^{n}\log lik(\theta,\widehat\eta(\theta))(X_i)$. Define
\begin{eqnarray}
\widehat\theta_n=\arg\sup_{\theta\in\Theta}\widehat S_n(\theta).\label{semiest}
\end{eqnarray}
In view of the above discussions, we can show that $\widehat\theta_n$ is semiparametric efficient if $\widehat\eta(\theta)$ is a consistent estimate of $\eta_\ast(\theta)$. The technical derivations in the above can be referred to Section 4 of \cite{sw92}. However, the form of $\widehat\theta_n$ depends on how we estimate the abstract $\eta(\theta)$ defined in (\ref{lfcurve}). For example, $\widehat\theta_n$ is just the semiparametric MLE if $\widehat\eta(\theta)$ is the well defined NPMLE. When the infinite dimensional $\mathcal{H}$ is too large, we may consider estimating $\eta_\ast(\theta)$ under some form of regularization, e.g., penalization. It is well known that the convergence rate of $\widehat\eta(\theta)$ is determined by the size of $\mathcal{H}$ in terms of its entropy number and the smoothing parameters associated with regularization methods (if used), e.g., smoothing parameter in penalized estimation.

In the following, we will consider two types of $\widehat\theta_n$ defined in (\ref{semiest}) according to how we estimate $\eta_\ast(\theta)$: (i) pure nonparametric estimation in Section~\ref{semisec}; (ii) nonparametric estimation under regularization in Section~\ref{smre}. Define $R_n\asymp r_n$ if $r_n/M\leq R_n\leq r_nM$ for some $M\geq 1$. We use $\mathcal{N}(\theta_0)$ to denote a neighborhood of $\theta_0$. Let $v_{i}$ denote the $i$-th unit vector in $\mathbb{R}^{d}$. Define the $i$-th ($(i,j)$-th) element of a vector $V$ $(\mbox{Matrix}\;M)$ as $V_{i}$ $(M_{ij})$. For a tensor $T^{(3)}(\theta)$, we define $V^{T}\otimes T^{(3)}(\theta)\otimes V$ as a $d$-dimensional vector with $i$-th element $V^{T}(\partial^{2}/\partial \theta^{2})(\dot T(\theta))_{i}V$, where $\dot T(\theta)$ is the first derivative of $T(\theta)$. Denote $int[x]$ and $\widetilde{int}[x]$ as the smallest nonnegative integer $\geq x$ and $>x$, respectively. The symbols $\mathbb{P}_n$ and $\mathbb{G}_n\equiv\sqrt{n}(\mathbb{P}_n-P)$ are used
for the empirical distribution and the empirical process of the observations,
respectively.

\section{Semiparametric Maximum Likelihood Estimation}\label{semisec}
In this section, we consider the maximum likelihood estimation of $\theta$ which corresponds to the case that (i) $\widehat\eta(\theta)$ is the NPMLE for $\eta_\ast(\theta)$ given any fixed $\theta$ and (ii) $\widehat S_n(\theta)=\log pl_n(\theta)$. The pure nonparametric estimation of $\eta_\ast(\theta)$ is often feasible when $\eta$ is under shape restrictions, e.g. the monotone cumulative hazard function. In general, the profile likelihood does not have a closed form since it is defined as a supremum over an infinite dimensional parameter space, see (\ref{prolik}). In practice, it can only be calculated numerically, e.g., via the iterative convex minorant algorithm \cite{h96}. We first discuss the construction of $\widehat\theta_n^{(k)}$, and then show that the minimal number of iterations $k^\ast$ is jointly determined by the convergence rates of $\widehat\theta_n^{(0)}$ and $\widehat\eta(\theta)$. In the end, two classes of semiparametric models are presented to illustrate our theories.

Throughout this section, we assume the following convergence rate Condition (\ref{convrate}) and the LFS Conditions M1-M4 specified in Appendix. For any random sequence $\widetilde{\theta}_n\overset{P}{\rightarrow}\theta_0$, we assume that
\begin{eqnarray}
\|\widehat{\eta}(\widetilde{\theta}_{n})-\eta_{0}\|=
O_{P}(\|\widetilde{\theta}_{n}-\theta_{0}\|\vee
n^{-r}),\label{convrate}
\end{eqnarray}
where $\|\cdot\|$ is some norm in $\mathcal{H}$ and $1/4<r\leq 1/2$. Of course we take the largest such $r$ in the following and call it the convergence rate for estimating $\eta$. The above range of $r$ holds in regular semiparametric models, which we can define without loss of generality to be models where the entropy integral converges. Theorems 3.1-3.2 in \citep{mv99} can be applied to calculate the convergence rate (\ref{convrate}). Under the above regularity conditions, Cheng and Kosorok (2008b) showed the following second order asymptotic linear expansion result.
\begin{theorem}\label{asynorlem}
Suppose that Conditions M1-M4 and (\ref{convrate})
hold. Also suppose that the MLE $\widehat{\theta}_n$ is consistent and $\widetilde I_0$ is nonsingular. We have
\begin{eqnarray}
\sqrt{n}(\widehat{\theta}_{n}-\theta_{0})=
\frac{1}{\sqrt{n}}\sum_{i=1}^{n}\widetilde{I}_{0}^{-1}\widetilde{\ell}_{0}
(X_{i})+O_P(n^{-2r+1/2}).\label{mleexp}
\end{eqnarray}
\end{theorem}

We need to estimate $\mathbb{P}_n\widetilde\ell_0$ and $\widetilde I_0$ to construct $\widehat\theta_n^{(k)}$ generated from the NR algorithm. In view of (\ref{lfs}) and (\ref{gpl-info2}), we can estimate them based on the derivatives of the log-profile likelihood (the sample analog of $S_n(\theta)$) as follows
\begin{eqnarray}
\left[\widehat{\ell}_{n}(\theta,s_{n})\right]_{i}&=&\frac{\log
pl_{n}(\theta+s_{n}v_{i})-\log
pl_{n}(\theta)}{ns_{n}},\label{estesco}\\
\left[\widehat{I}_{n}(\theta,t_{n})\right]_{i,j}&=&-\frac{\log
pl_{n}(\theta+t_{n}(v_i+v_{j}))+\log
pl_{n}(\theta)}{nt_{n}^{2}}\nonumber\\&&+\frac{\log
pl_{n}(\theta+t_{n}v_i)+\log
pl_{n}(\theta+t_{n}v_j)}{nt_{n}^{2}}.\label{estei}
\end{eqnarray}
In the above we use the numerical derivatives since the smoothness and differentiability of $\log pl_n(\theta)$ are usually unknown. In Lemma~\ref{mainthm} of Appendix, we show that (\ref{estesco}) and (\ref{estei}) (also called as the observed information in \citep{mv99}) are indeed the consistent estimators. Thus, we can write $\widehat\theta_n^{(k)}$ in step (III) as
\begin{eqnarray}
\widehat{\theta}_{n}^{(k)}=\widehat{\theta}_{n}^{(k-1)}+
\left[\widehat I_{n}\left(\widehat{\theta}_{n}^{(k-1)},t_n^{(k-1)}\right)\right]^{-1}
\widehat\ell_{n}\left(\widehat{\theta}_{n}^{(k-1)},s_n^{(k-1)}\right), \label{sche}
\end{eqnarray}
where step sizes $s_n^{(k-1)}\vee t_n^{(k-1)}=o(1)$. A close inspection of (\ref{sche}) reveals that we have constructed $\widehat\theta_n^{(k)}$ even without knowing the forms of $\widetilde\ell_0$ and $\widetilde I_0$. 

The convergence of $\widehat\theta_n^{(k)}$ to $\widehat\theta_n$, which is exactly the maximizer of $\log pl_n(\theta)$, as $k\rightarrow\infty$ is guaranteed by the asymptotic parabolic form of $\log pl_n(\theta)$ proven in \cite{mv00}. However, to figure out the minimal $k^\ast$ such that $\|\widehat\theta_n^{(k^\ast)}-\widehat\theta_n\|=o_P(n^{-1/2})$, we need to make use of the second order asymptotic quadratic expansion of $\log pl_n(\theta)$ derived in \cite{ck08b} under the above regularity conditions. As seen from (\ref{sche}), the orders of step sizes $(s_n^{(k-1)},t_n^{(k-1)})$ are critical in determining the convergence rate of $\widehat\theta_n^{(k)}$ to $\widehat\theta_n$, and thus need to be properly chosen at each iteration. In the below Lemma, we present the optimal step sizes, under which the fastest convergence rate is achieved, at each iteration. Denote the convergence rate of $\|\widehat\theta_n^{(k-1)}-\widehat\theta_n\|$ as $O_P(n^{-r_{k-1}})$.
\begin{lemma}\label{step}
Suppose Conditions in Theorem~\ref{asynorlem} hold. The convergence rate of $\|\widehat\theta_n^{(k)}-\widehat\theta_n\|$ is improved through the following three stages: 
\begin{itemize}
\item[(i)] $\|\widehat\theta_n^{(k)}-\widehat\theta_n\|=O_P(\|\widehat\theta_n^{(k-1)}-\widehat\theta_n\|^{3/2})$ when $r_{k-1}<r$ and we choose $(s_n^{(k-1)},t_n^{(k-1)})\asymp(n^{-3r_{k-1}/2},n^{-r_{k-1}/2})$;
    
\item[(ii)] $\|\widehat\theta_n^{(k)}-\widehat\theta_n\|=O_P(\|\widehat\theta_n^{(k-1)}-\widehat\theta_n\|^{1/2}n^{-r})$ when $r\leq r_{k-1}< 1/2$ and we choose $(s_n^{(k-1)},t_n^{(k-1)})\asymp(n^{-r-r_{k-1}/2}, n^{-r_{k-1}/2})$;
    
\item[(iii)] $\|\widehat\theta_n^{(k)}-\widehat\theta_n\|=O_P(n^{-r-1/4})$ when $r_{k-1}\geq 1/2$ and we choose $(s_n^{(k-1)},t_n^{(k-1)})\asymp(n^{-r-1/4}, n^{-r_{k-1}/2})$.
\end{itemize}
\end{lemma}
Now we present our first main theorem, i.e., Theorem~\ref{thm1}. Let $\widehat\theta_n^{(0)}$ be $n^\psi$-consistent. We first show that $\|\widehat{\theta}_n^{(k)}-\widehat{\theta}_n\|=O_P(n^{-S(\psi,r,k)})$ based on which we figure out the value of $k^\ast$  in (\ref{kstar}).  According to the above Lemma~\ref{step}, it is easily seen that $S(1/2,r,k)=r+1/4$ for any $1/4<r\leq 1/2$ and $k\geq 1$ (thus $k^\ast=1$); and $S(1/3,1/2,1)=1/2$ and $S(1/3,1/2,k)=3/4$ for any $k\geq 2$ (thus $k^\ast=2$). Following similar logic, we can give the general form of $S(\psi,r,k)$ as follows. Define, if $\widetilde S_1(\psi,r)\geq 1/2$,
\begin{eqnarray*}
S(\psi,r,k)=\left\{
\begin{array}{lr}
S_1(\psi,k)&k\leq K_1(\psi,r)\\
r+1/4&k\geq K_1(\psi,r)+1
\end{array}
\right.,\label{sform0}
\end{eqnarray*}
where $S_1(\psi,k)=\psi(3/2)^k$, $K_1(\psi,r)=int\left[\log(r/\psi)/\log(3/2)\right]$ and $\widetilde S_1(\psi,r)=S_1(\psi,K_1(\psi,r))$, and if 
 $r\leq\widetilde S_1(\psi,r)<1/2$, 
\begin{eqnarray*}
S(\psi,r,k)=\left\{
\begin{array}{lr}
S_1(\psi,k)&k\leq K_1(\psi,r)\\
S_2(\widetilde S_1(\psi,r),r,k-K_1(\psi,r))& K_1(\psi,r)<k\leq K_1(\psi,r)+\widetilde K_2(\psi,r)\\
r+1/4&k\geq K_1(\psi,r)+\widetilde K_2(\psi,r)+1
\end{array}
\right.,\label{sform}
\end{eqnarray*}
where $S_2(\psi,r,k)=
2r+2^{-k}(\psi-2r)$, $K_2(\psi,r)=int
[\log\{(2r-\psi)/(2r-1/2)\}/\log
2]$ and $\widetilde K_2(\psi,r)=K_2(\widetilde S_1(\psi,r),r)$.

\begin{theorem}\label{thm1}
Suppose that Conditions in Theorem~\ref{asynorlem} hold and proper step sizes are chosen according to Lemma~\ref{step}. Let $\widehat{\theta}_{n}^{(k)}$ be the k-step estimator defined in (\ref{sche}) and $\widehat{\theta}_{n}^{(0)}$ be $n^{\psi}$-consistent for
$0<\psi\leq 1/2$. Recall that $\|\widehat{\theta}_n^{(k)}-\widehat{\theta}_n\|=O_P(n^{-r_k})$. We show that $r_k$ increases from $\psi$ to $(r+1/4)$ as $k\rightarrow\infty$. Specifically, we have
\begin{eqnarray}
\|\widehat{\theta}_n^{(k)}-\widehat{\theta}_n\|
=O_P(n^{-S(\psi,r,k)}).\label{form}
\end{eqnarray}
This implies that
\begin{eqnarray}
\|\widehat\theta_n^{(k^\ast)}-\widehat\theta_n\|=o_P(n^{-1/2}),\label{kstar}
\end{eqnarray}
where $k^\ast=K_1(\psi,r)+\widetilde{int}[\log((2r-\widetilde S_1(\psi,r))/(2r-1/2))/\log 2]$.
\end{theorem}
Interestingly, we notice that the optimal bound of $\|\widehat{\theta}_n^{(k)}-\widehat{\theta}_n\|$, i.e. $O_P(n^{-r-1/4})$, is intrinsically determined by how accurately we estimate the nuisance parameter, i.e., the value of $r$. This bound can not be further improved unless we are willing to make stronger
assumptions than M1-M4, which seem unrealistic. From the form of $S(\psi,r,k)$, we find that more accurate initial estimate leads to higher order asymptotic efficiency of $\widehat{\theta}_n^{(k)}$. How to obtain $\widehat\theta_n^{(0)}$ with proper convergence rate will be discussed in Section~\ref{iniest}.

We apply Theorem~\ref{thm1} to the following two examples whose detailed technical illustrations and model assumptions can be found in \cite{mv00,rcl96}. The required Conditions in Theorem~\ref{thm1} are verified in \cite{ck08a,ck08b} for Examples 1-2. We can also apply our theory to the semiparametric regression model under shape constraints, e.g., \cite{c09}.

{\it Example 1: Cox Model under Current Status Data}

In the Cox proportional hazards model, the hazard function of the
survival time $T$ of a subject with covariate $Z$ is expressed
as:
\begin{eqnarray*}
\lambda(t|z)\equiv\lim_{\Delta\rightarrow 0}\frac{1}{\Delta}Pr(t\leq T<t+\Delta|T\geq t,Z=z)=\lambda(t)\exp(\theta' z),
\end{eqnarray*}
where $\lambda$ is an unspecified baseline hazard function. We consider the current status data where each subject is observed at a single examination time $Y$ to determine if an event has occurred, but the event time $T$ cannot be known exactly. Specifically, the observed data are $n$ realizations of
$X=(Y, \delta, Z)\in R^{+}\times \lbrace 0,1 \rbrace\times R$, where $\delta=I\{T \leq Y\}$. The cumulative hazard function
$\eta(y)=\int_{0}^{y}\lambda(t)dt$ is considered as the nuisance parameter. The parameter space $\mathcal{H}$ for $\eta$ is
restricted to a set of nondecreasing and cadlag functions on some compact interval. In this model, it is well known that both $\widehat\eta(\theta)$ and $\log pl_n(\theta)$ have no explicit forms, and can only be calculated numerically via the iterative convex minorant algorithm, see \cite{h96}. As for the convergence rate of $\eta$, Murphy and van der Vaart (1999) showed
$\|\widehat{\eta}(\widetilde{\theta}_{n})-\eta_{0}\|_{2}=O_{P}(\|\widetilde{\theta}_{n}-\theta_{0}\|\vee n^{-1/3})$, where $\|\cdot\|_{2}$ is the $L_{2}$ norm. According to Theorem~\ref{thm1}, we establish the following table to depict the convergence of $\widehat\theta_n^{(k)}$ to $\widehat\theta_n$ given different initial estimates until it reaches the lower bound $O_P(n^{-7/12})$.
\begin{center}
Table 1. {\it Cox Model under Current Status Data $(r=1/3)$ }\vspace{0.1in} 
\centering 
\begin{tabular}{c c c c}
\hline\hline
 & $\psi=1/2$ & $\psi=1/3$ & $\psi=1/4$\\
\hline
\mbox{Cox} &  $r_1=7/12$ & $r_1=1/2, r_2=7/12$ & $r_1=3/8, r_2=25/48, r_3=7/12$\\
\mbox{Models}& $k^\ast=1$ & $k^\ast=2$ & $k^\ast=2$\\
\multicolumn{4}{@{}p{12.6cm}@{}}{\rule{12.9cm}{0.2pt}}\\
\multicolumn{4}{@{}p{12.6cm}@{}}
{\small Remark: Define $\|\widehat\theta_n^{(k)}-\widehat\theta_n\|=O_P(n^{-r_k})$.}
\end{tabular}
\end{center}

{\it Example 2: Semiparametric Mixture Model in Case-Control Studies}

Roeder, Carroll and Lindsay (1996) consider the logistic regression model with a missing covariate for case-control studies. In this model, they observe two independent random samples: one complete component $Y_{C}=(D_{C},W_{C})$ and $Z_C$ of the size $n_{C}$, and one reduced component $Y_{R}=(D_{R},W_{R})$ of the size $n_{R}$. Following the assumptions given in \cite{rcl96}, the
likelihood for $x=(y_C,y_R,z_C)$ is defined as
\begin{eqnarray*}
lik(\theta',\eta)(x)=p_{\theta'}(y_{C}|z_{C})\eta\{z_{C}\}\int
p_{\theta'}(y_{R}|z)d\eta(z),
\end{eqnarray*}
where $d\eta$ denotes the density of $\eta$ w.r.t. some dominating measure, and
\begin{eqnarray*}
p_{\theta'}(y|z)=\left(\frac{\exp(\gamma+\theta e^{z})}{1+\exp(\gamma+\theta e^{z})}\right)^{d}\left(\frac{1}{1+\exp(\gamma+\theta e^{z})}\right)^{1-d}\phi_\sigma(w-\alpha_{0}
-\alpha_{1}z),
\end{eqnarray*}
where $\phi_{\sigma}(\cdot)$ denotes the density for $N(0,\sigma)$. The unknown parameters are
$\theta'=(\theta,\alpha_{0},\alpha_{1},\gamma,\sigma)$ ranging over the
compact $\Theta'\subset\mathbb{R}^{4}\times(0,\infty)$ and the
distribution $\eta$ of the regression variable restricted
to the set of nondegenerate probability distributions with a known compact support. In this semiparametric mixture model, we will concentrate on the
regression coefficient $\theta$, considering
$\theta_{2}=(\alpha_{0},\alpha_{1},$ $\gamma,\sigma)$
and $\eta$ as nuisance parameters. The NPMLE $\widehat\eta(\theta)(z)$ is a weighted average of two empirical distributions and
the log-profile likelihood implicitly defined as
$$\widehat S_n(\theta)=\log pl_n(\theta)=\sup_{\theta_2,\eta}\log lik_n(\theta',\eta)$$ has no explicit form. Let $(\widehat{\theta}_{2,\theta},\widehat{\eta}(\theta))$ be the profile
likelihood estimator for $(\theta_{2},\eta)$
so that $\widehat{\theta}_{\theta}'=(\theta,\widehat{\theta}_{2,\theta})$. Both $\widehat\eta(\theta)$ and $\widehat S_n(\theta)$ can be computed efficiently via the iterative algorithm in Section 4 of \cite{rcl96}, a special case of our general algorithm.
Murphy and van der Vaart (1999) showed that, for any $\widetilde{\theta}_{n}\overset{P}{\rightarrow}\theta_{0}$,
\begin{eqnarray}
\|\widehat{\eta}(\widetilde{\theta}_{n})-\eta_{0}\|_{BL_{1}}+\|\widehat{\theta}_{\widetilde{\theta}_{n}}'-\theta_{0}'\|=
O_P(|\widetilde{\theta}_{n}-\theta_{0}|\vee n^{-\frac{1}{2}}),\label{eg3rate}
\end{eqnarray}
where $\|\cdot\|_{BL_1}$ is the weak topology. This implies that $r=1/2$. The following Table 2 is similar as Table 1. Interestingly, we find that $\widehat\theta_n^{(k)}$ converges to $\widehat\theta_n$ at a faster rate in the second model.
\begin{center}
Table 2. {\it Semiparametric Mixture Model in Case-Control Studies $(r=1/2)$}\vspace{0.1in} 
\centering 
\begin{tabular}{c c c c}
\hline\hline
 & $\psi=1/2$ & $\psi=1/3$ & $\psi=1/4$\\
\hline
\mbox{Mixture} & $r_1=3/4$ & $r_1=1/2, r_2=3/4$ & $r_1=3/8, r_2=9/16, r_3=3/4$\\
\mbox{Models}& $k^\ast=1$ & $k^\ast=2$ & $k^\ast=2$\\
\multicolumn{4}{@{}p{12.0cm}@{}}{\rule{12.3cm}{0.2pt}}\\
\multicolumn{4}{@{}p{12.0cm}@{}}
{\small Remark: Define $\|\widehat\theta_n^{(k)}-\widehat\theta_n\|=O_P(n^{-r_k})$.}
\end{tabular}
\end{center}
\section{Semiparametric Estimation under Regularization}\label{smre}
In this section, we consider the semiparametric estimation under two types of regularizations, i.e., kernel estimation and penalized estimation. In contrast with the profile likelihood estimation, the regularized $\widehat S_n(\theta)$ is usually differentiable although its form may vary under different regularizations. We first present a unified framework for studying $\widehat \theta_n^{(k)}$ when $\widehat S_n(\theta)$ is third order differentiable, and then present several examples corresponding to different regularizations which fit into this framework. We also discuss the variable selection in partly linear models as an extension of the penalized estimation.

In this section, we construct $\widehat\theta_n^{(k)}$ in step (III) as follows:
\begin{eqnarray}
\widehat{\theta}_n^{(k)}=\widehat{\theta}_n^{(k-1)}+\left[\widehat
I_n(\widehat{\theta}_n^{(k-1)})\right]^{-1}\widehat\ell_n(\widehat{\theta}_n^{(k-1)}),\label{parasche-0}
\end{eqnarray}
where $\widehat\ell_n(\cdot)=\widehat S_n^{(1)}(\cdot)/n$ and
\begin{eqnarray}
\widehat I_n(\cdot)=-\widehat S_n^{(2)}(\cdot)/n,\label{semicon1}
\end{eqnarray}
where $\widehat S_n^{(j)}(\cdot)$ is the $j$-th derivative of $\widehat S_n(\cdot)$. When $\widehat S_n^{(2)}(\theta)$ has no explicit form or is hard to compute, we may prefer constructing $[\widehat I_{n}(\theta)]_{ij}$ as
\begin{eqnarray}
-n^{-1/2}\frac{
[\widehat{S}_n^{(1)}(\theta+n^{-1/2}t_2v_j)]_i-[\widehat{S}_n^{(1)}
(\theta+n^{-1/2}t_1v_j)]_i}{t_2-t_1},\label{semicon2}
\end{eqnarray}
where $t_1$ and $t_2$ $(t_1<t_2)$ are arbitrarily fixed real numbers.

Recall that $$S_n(\theta)=n\mathbb{P}_n\log lik(\theta,\eta_\ast(\theta))$$ and define $S_n^{(j)}(\cdot)$ as the $j$-th derivative of $S_n(\cdot)$. In view of the discussions in Section~\ref{secrev}, i.e. (\ref{lfs}) \& (\ref{gpl-info2}), we expect that $\widehat\theta_n^{(k)}$ converges to $\widehat\theta_n$ if $\widehat S_n^{(j)}(\cdot)$ approximates $S_n^{(j)}(\cdot)$ well enough round $\theta_0$ for $j=1,2,3$. Therefore, we assume the following general condition G.
\begin{enumerate}
\item[G.] Assume that
\begin{eqnarray}
\frac{1}{n}\widehat{S}^{(1)}_n(\theta_0)-\frac{1}{n}
S_n^{(1)}(\theta_0)&=&O_P(n^{-2g}),\label{pricon2}\\
\sup_{\theta\in\mathcal{N}(\theta_0)}\left|\frac{1}{n}\widehat{S}^{(2)}_n
(\theta)-\frac{1}{n} S_n^{(2)}(\theta)\right|&=&
O_P(n^{-g}),\label{pricon3}\\
\sup_{\theta\in\mathcal{N}(\theta_0)}\left|\frac{1}{n}\widehat
S_n^{(3)}(\theta)\right|
&=&O_P(1),\label{pricon4}
\end{eqnarray}
where $1/4<g\leq 1/2$.
\end{enumerate}
We will provide two sets of sufficient conditions for G in the kernel estimation, where the value of $g$ is determined by the bandwidth order,
and in the penalized estimation, where the value of $g$ is determined by the smoothing parameter, respectively. In this sense, we can think $g$ is a measure of the convergence rate of $\widehat\eta(\theta)$ as in (\ref{convrate}). We may verify (\ref{pricon4}) by showing
\begin{eqnarray}
\sup_{\theta\in\mathcal{N}(\theta_0)}\left|\frac{1}{n}\widehat
S_n^{(3)}(\theta)-\frac{1}{n}S_n^{(3)}(\theta)\right|
&=&o_P(1),\label{pricon5}
\end{eqnarray}
and that the class of functions
$\{(\partial^3/\partial\theta^3)\log
lik(x;\theta,\eta_\ast(\theta)):\theta\in\mathcal{N}(\theta_0)\}$ is
P-Glivenko-Cantelli and that
$$\sup_{\theta\in\mathcal{N}(\theta_0)}E\left|(\partial^3/\partial\theta^3)\log
lik(X;\theta,\eta_\ast(\theta))\right|<\infty.$$

Now we present our second main theorem, i.e., Theorem~\ref{conv-gpl-kmle}. Define
\begin{align}\label{rform}
R(\psi,g,k)=
\begin{cases}
R_1(\psi,g,k)& \text{$k\leq L_1(\psi, g)$}\\
R_2(R_1(\psi,g,L_1(\psi,g)),g,k-L_1(\psi,g))& \text{
$k>L_1(\psi, g)$}
\end{cases}
\end{align}
where $R_1(\psi,g,k)=(1/2-g)+2^k(\psi+g-1/2)$, $L_1(\psi,g)=int[\log (g/(g+\psi-1/2))/\log
2]$, $\widetilde L_1(\psi,g)=\widetilde{int}[\log(g/(g+\psi-1/2))/\log
2]$ and $R_2(\psi,g,k)=kg+\psi$.
\begin{theorem}\label{conv-gpl-kmle}
Suppose that Condition G holds, $\widehat{\theta}_{n}$ defined in (\ref{semiest}) is consistent and $\widetilde{I}_0$ is nonsingular. We have
\begin{eqnarray}
\sqrt{n}(\widehat{\theta}_{n}-\theta_{0})=
\frac{1}{\sqrt{n}}\sum_{i=1}^{n}\widetilde{I}_{0}^{-1}\widetilde{\ell}_{0}
(X_{i})+O_P(n^{1/2-2g}).\label{gpl-mleexp}
\end{eqnarray}
Let
$\widehat{\theta}_{n}^{(k)}$ be the $k$-step estimator defined in
(\ref{parasche-0}) and $\widehat{\theta}_{n}^{(0)}$ be $n^\psi$-consistent for
$(1/2-g)<\psi\leq 1/2$. Define $\|\widehat{\theta}_{n}^{(k)}-\widehat{\theta}_{n}\|=O_P(n^{-r_k})$.
We show that $r_k$ increases from $\psi$ to $\infty$ as $k\rightarrow\infty$. Specifically, we show
\begin{eqnarray}
\|\widehat{\theta}_{n}^{(k)}-\widehat\theta_{n}\|&=& O_P(n^{-2^k
\psi})\;\;\;\;\;\;\;\;\mbox{if}\;\;\widehat{I}_{n}(\cdot)\; \mbox{is
defined in (\ref{semicon1})}
\label{gpl-thm-1},\\
\|\widehat{\theta}_{n}^{(k)}-\widehat\theta_{n}\|&=&
O_P(n^{-R(\psi,g,k)})\;\;\mbox{if}\;\;\widehat{I}_{n}(\cdot)\;
\mbox{is defined in (\ref{semicon2})}.\label{gpl-thm-2}
\end{eqnarray}
This implies that
$$\|\widehat\theta_n^{(k^\ast)}-\widehat\theta_n\|=o_P(n^{-1/2}),$$ where $k^\ast=\widetilde{int}[\log (1/2\psi)/\log 2]$ for (\ref{gpl-thm-1}) and  $k^\ast=\widetilde{L}_1(\psi,g)$ for (\ref{gpl-thm-2}).
\end{theorem}
Note that (\ref{gpl-thm-1}) is a
statistical counterpart to the well known quadratic convergence of the Newton-Raphson algorithm; see Page 312 of \citep{or70}. Theorems~\ref{thm1} and~\ref{conv-gpl-kmle} imply that (i) more than $k^\ast$ iterations, i.e., $k$, will not change the limiting distribution of $\widehat\theta_n^{(k)}$, but will improve its higher order asymptotic efficiency; (ii) the higher order asymptotic efficiency of $\widehat\theta_n^{(k)}$ is determined by how accurately $\eta$ is estimated, i.e., the values of $r$ and $g$; (iii) $\widehat\theta_n^{(k)}$ converges to $\widehat\theta_n$ faster when $\widehat I_n$ is constructed as an analytical derivative no matter whether the regularization is used or not.

\begin{remark}
Given that the initial estimate is $\sqrt{n}$ consistent, we have
\begin{eqnarray*}
\|\widehat{\theta}_{n}^{(k)}-\widehat{\theta}_{n}\|&=&
O_P(n^{-2^{k-1}})\;\;\;\;\;\;\;\mbox{if}\;\;\widehat{I}_{n}(\cdot)\;
\mbox{is constructed
as in (\ref{semicon1}),}\nonumber\\
\|\widehat{\theta}_{n}^{(k)}-\widehat{\theta}_{n}\|&=&
O_P(n^{-(1/2+kg)})\;\;\mbox{if}\;\;\widehat{I}_{n}(\cdot)\; \mbox{is
constructed as in (\ref{semicon2})}
\end{eqnarray*}
based on Theorem~\ref{conv-gpl-kmle}. This implies $k^\ast=1$.
\end{remark}

A by-product of Theorem~\ref{conv-gpl-kmle} is the application to the parametric models, i.e., $\eta$ is known. In this case, $\widehat S_n(\theta)$ becomes $\ell_\theta(X)=\log lik(\theta; X)$, and we simplify the general Condition G to the following Conditions P1-P2. Denote the first, second and third derivative of
$\ell_\theta(\cdot)$ w.r.t. $\theta$ as
$\dot{\ell}_\theta(\cdot)$, $\ddot{\ell}_\theta(\cdot)$ and
$\ell^{(3)}_\theta(\cdot)$, respectively. The information matrix at $\theta_0$ is defined as $I_0$.
\begin{enumerate}
\item[P1.] $\dot{\ell}_\theta(\cdot)$ and
$\ddot{\ell}_\theta(\cdot)$ are absolutely continuous in $\theta$.

\item[P2.] There exists a $\delta>0$ such that, for any
$|t|\leq\delta$,
\begin{eqnarray}
E\left[\ell_{\theta_0+t}^{(i+1)}(X_1)\right]^2\leq
K_i\;\;\;\mbox{for some finite constant}\;K_i,\label{cond2}
\end{eqnarray}
where $i=1,2$.
\end{enumerate}
We can easily prove Corollary~\ref{conv-para} by following similar analysis in Theorem~\ref{conv-gpl-kmle} and considering Lemma~\ref{lemm-para}. Thus, its proof is skipped.

\begin{corollary}\label{conv-para}
Suppose that Conditions P1 \& P2 hold. Also suppose that the parametric MLE $\widehat{\theta}_n$ is consistent and $I_0$ is nonsingular. Let $g=1/2$. Then all the conclusions for $\widehat\theta_n$ and $\widehat\theta_n^{(k)}$ in Theorem~\ref{conv-gpl-kmle} hold for the parametric estimation.
\end{corollary}
The above corollary generalizes the one/two-step parametric estimation results in \cite{jjv85}.
Comparing Theorem~\ref{conv-gpl-kmle} with Corollary~\ref{conv-para}, we notice that $\widehat\theta_n^{(k)}$ converges to $\widehat\theta_n$ at a slower rate in semiparametric models. This results from the presence of an infinite dimensional $\eta$ estimated at a slower-than-parametric rate by comparing Lemmas~\ref{gpl-le3} and~\ref{lemm-para}.

\begin{remark}\label{remak1}
We would like to mention that the regularized $\widehat S_n(\theta)$ may not be differentiable in some semiparametric models, e.g., the penalized estimation of partly linear models under current status data studied in \cite{ck09}. In such cases, we can take the discretization approach to construct $\widehat\theta_n^{(k)}$ as in the profile likelihood framework, i.e., (\ref{sche}), and obtain similar results as in Theorem~\ref{thm1} if we can prove that the non-smooth $\widehat S_n(\theta)$ share the same higher order quadratic expansion as $\log pl_n(\theta)$. Indeed, Cheng and Kosorok (2009) have proven such results for the non-smooth regularized $\widehat S_n(\theta)$ under weaker conditions. See \cite{c10} for more elaborations.
\end{remark}

\subsection{Kernel Estimation in Semiparametric Models}\label{kerest0}
In this subsection, we consider the kernel estimation in semiparametric models. Due to its simple form, the kernel estimate of $\eta$ and the related iterative algorithm of estimating $\theta$ are widely used in semiparametric models, e.g., \cite{a95,s88}. In particular, the kernel approach is proven to be a powerful inferential tool for the class of conditionally parametric models (CPM), see \cite{sw92, ss94}. Thus, in this subsection, we will focus on the class of CPM although our conclusions can be extended to more general class of semiparametric models by incorporating the results in \cite{a95}. Under kernel estimation, $k^\ast$ is shown to depend on the order of bandwidth used in the kernel function.

The class of CPM was first introduced by Severini and Wong (1992) and further generalized to the quasi-likelihood framework by Severini and Staniswalis (1994). Specifically, we observe $X=(Y,W,Z)$ such that the distribution of $Y$ conditional on partitioned covariates $W=w$ and $Z=z$ is parameterized by a finite dimensional parameter $\phi=(\theta, \lambda_z)$, where $\lambda_z\in H\subset\mathbb{R}$ depends on the value of $z$ as a function $\eta(z)$. The joint distribution of $(W,Z)$ is assumed to be independent of $\phi$. Thus, this semiparametric model  has the log-likelihood $\log lik(X;\theta,\eta(z))$ and is called conditionally parametric. The practical performance of the iterative estimation procedure (I)-(IV) for the CPM is extensively studied in \cite{ss94}.

We assume that $\eta(z)\in\mathcal{H}=\{h\in C^2(\mathcal{Z}): h(z)\in \mbox{interior}(H)\;\mbox{for all}
\;z\in\mathcal{Z}\}$. An important feature of CPM is that its least favorable curve can be expressed as (see \cite{sw92} for details)
\begin{eqnarray}
\eta_\ast(\theta)(z)=\arg\sup_{\eta\in
C^2[0,1]}E[\log lik(X;\theta,\eta)|Z=z],\label{kerest1}
\end{eqnarray}
and thus
its kernel estimate is written as
\begin{eqnarray}
\widehat{\eta}(\theta)(z)=\arg\sup_{\eta\in
C^2[0,1]}\sum_{i=1}^{n}\log lik(X_i;\theta,\eta(Z_i))
K\left(\frac{z-Z_{i}}{b_n}\right),\label{kerest}
\end{eqnarray}
where $K(\cdot)$ is a kernel with the bandwidth $b_n\rightarrow
0$. For example, if $(Y|w=W,Z=z)\sim N(\theta'w,\eta(z))$, then we have
\begin{eqnarray}
\widehat{\eta}(\theta)(z)=\frac{\sum_{i=1}^{n}(Y_i-\theta'W_i)^2
K((z-Z_i)/b_n)}{\sum_{i=1}^{n}K((z-Z_i)/b_n)},\nonumber\\
\widehat{S}_n(\theta)=-\frac{1}{2}\sum_{i=1}^{n}
\frac{(Y_i-\theta'W_i)}{\widehat{\eta}(\theta)(Z_i)}-\frac{1}{2}
\sum_{i=1}^{n}\log \widehat{\eta}(\theta)(Z_i).\label{eg1gpl}
\end{eqnarray}
Although $\widehat\eta(\theta)$ (and thus $\widehat S_n(\theta)$) solved from (\ref{kerest}) generally has no explicit form, based on (\ref{kerest}) we can control the asymptotic behaviors of $\widehat\eta(\theta)$ (and thus $\widehat S_n(\theta)$) by assuming proper kernel conditions, see the below Example 3.

By exploiting the parametric structure of CPM, we will show $\widehat S_n(\theta)$ satisfies the general Condition G under the below Conditions K1-K2 and C1-C2.
\begin{enumerate}
\item[K1.] For arbitrary $\theta_1\in\Theta$ and $\lambda_1\in H$,
if $\theta\neq\theta_1$, then $E_{\theta_1,\lambda_1}\log
lik(X;\theta,\lambda)<E_{\theta_1,\lambda_1} \log
lik(X;\theta_1,\lambda_1)$;

\item[K2.] Assume that
\begin{eqnarray}
E\left\{\sup_{(\theta,\lambda)\in\Theta\times
H}\left|\frac{\partial^{r+s}\log
lik(X;\theta,\lambda)}{\partial\theta^r
\partial\lambda^s}\right|^2\right\}<\infty
\end{eqnarray}
for all $r,s=0,\ldots,4$ and $r+s\leq 4$.
\end{enumerate}
Similar identifiability Condition K1 and smoothness Condition K2 are also used in \cite{sw92}. Our next conditions C1-C2 are concerned about the smoothness and convergence rate of $\eta_\ast(\theta)$ and $\widehat\eta(\theta)$.
We denote the derivative of $\eta_\ast(\theta)$ $(\widehat{\eta}(\theta))$ w.r.t. $\theta$ as $\eta_{\ast}^{(s)}(\theta)$ $(\widehat{\eta}^{(s)}(\theta))$, and their values at $\theta_0$ as $\eta_{\ast0}^{(s)}$ $(\widehat\eta_0^{(s)})$.
\begin{enumerate}
\item[C1.] Assume that, for all $r,s=0,1,2,3$ and $r+s\leq 3$,
$$\frac{\partial^{r+s}}{\partial z^r\partial\theta^s}
\eta_\ast(\theta)(z)\;\;\mbox{and}\;\; \frac{\partial^{r+s}}{\partial
z^r\partial\theta^s}\widehat{\eta}(\theta)(z)$$ exist and
$\sup_{\theta\in\mathcal{N}(\theta_0)}\|\eta_{\ast}^{(s)}(\theta)\|
_\infty<\infty$.
\end{enumerate}
\begin{enumerate}
\item[C2.] Assume that
\begin{eqnarray}
\sup_{\theta\in\mathcal{N}(\theta_0)}\|\widehat{\eta}^{(s)}(\theta)-
\eta^{(s)}_{\ast}(\theta)\|_\infty&=&O_P(n^{-g})
\;\;\;\;\;\;\mbox{for}\;s=0,1,2,\label{pricon10}\\
\sup_{\theta\in\mathcal{N}(\theta_0)}\|\widehat{\eta}^{(3)}(\theta)-
\eta^{(3)}_\ast(\theta)\|_\infty&=&o_P(1),\label{pricon50}\\
\left\|\frac{\partial}{\partial z}\widehat{\eta}_{0}(z)-
\frac{\partial}{\partial z}\eta_{\ast0}(z)\right\|_{\infty}&=&o_P(n^{-\delta}),\label{gplrate3}\\
\left\|\frac{\partial}{\partial z}\widehat{\eta}^{(1)}_{0}(z)-
\frac{\partial}{\partial
z}\eta^{(1)}_{\ast0}(z)\right\|_{\infty}&=&o_P(n^{-\delta}).\label{gplrate4}
\end{eqnarray}
for some $g\in(1/4,1/2]$ and $(2g-1/2)\leq \delta\leq g$.
\end{enumerate}

In view of (\ref{kerest1})-(\ref{kerest}), we can verify C2 by applying the kernel theories under some proper kernel conditions and K1-K2. For example, in Lemma~\ref{lem1}, we show that  the convergence rate of the kernel estimate in (\ref{pricon10}), which determines the value of $g$ in (\ref{pricon2})-(\ref{pricon3}), relies on the order of bandwidth $b_n$ used in (\ref{kerest}). Note that Condition C2 also satisfies (\ref{convrate}) assumed for the NPMLE since
\begin{eqnarray*}
\|\widehat{\eta}(\widetilde\theta_n)-\eta_0\|&\leq&
\|\widehat\eta(\widetilde\theta_n)-\widehat\eta(\theta_0)\|_{\infty}+\|\widehat\eta(\theta_0)-\eta_\ast(\theta_0)\|_\infty
\\&\leq&O_P(\|\widetilde\theta_n-\theta_0\|\vee n^{-g})
\end{eqnarray*}
by the construction that $\eta_\ast(\theta_0)=\eta_0$, C1-C2 and (\ref{pricon10}). Our conditions K1-K2 and C1-C2 are generally stronger than M1-M4 and (\ref{convrate}) since the semiparametric models under consideration have the assumed parametric structure.
\begin{theorem}\label{gpl-ver}
Assuming that Conditions K1-K2 and C1-C2 hold, then the
Condition G required in Theorem~\ref{conv-gpl-kmle} is
satisfied for the kernel estimation in conditionally parametric models.
\end{theorem}
The consistency of $\widehat\theta_n$ required in Theorem~\ref{conv-gpl-kmle} can be established if we further require the global condition $\sup_{\theta\in\Theta}\|\widehat{\eta}(\theta)-\eta_\ast(\theta)\|_{\infty}\rightarrow 0$, see Proposition 1 of \citep{sw92}. In the next example, we  apply Theorems~\ref{conv-gpl-kmle} and~\ref{gpl-ver} to a subclass of CPM, called conditionally exponential models (CEM), in which $\widehat\eta(\theta)$ has a closed-form. This makes the verifications of C1-C2 much easier. The relation between $k^\ast$ and the order of $b_n$ in (\ref{kerest}) is also specified in the below example. We may also apply our theories to the more complicated semiparametric transformation model, i.e., \cite{lsv08}.

{\it Example 3. Conditionally Exponential Models}

In CEM, there exists a function $\psi_\theta(\cdot)$ such that the conditional distribution of $\psi_\theta(Y,W)$ given $Z=z$ does not depend on $\theta$ and forms an exponential family. And its log-likelihood can be expressed as
$$\log lik(X;\theta,\eta)=\psi_\theta(Y,W)T(\eta(Z))-A(\eta(Z))+S(\psi_\theta(Y,W))$$ for some functions $T$, $A$ and $S$. Some simple algebra gives that
\begin{eqnarray}
\widehat\eta(\theta)(z)=\rho\left(\frac{\sum_{i=1}^{n}\psi_\theta(Y_i,W_i)
K((z-Z_i)/b_n)}{\sum_{i=1}^{n}K((z-Z_i)/b_n)}\right),\label{inter03}
\end{eqnarray}
where $\eta=\rho\{E_{\theta,\eta}(\psi_\theta(Y,W))\}$. In the previous conditional normal model, we have $\psi_\theta(Y,W)=(Y-\theta'W)^2$ and $\rho(t)=t$. Another example is that $(Y|W=w,Z=z)\sim\mbox{Exp}(0,\exp(\theta'w+\eta(z)))$ in which $\psi_\theta(Y,W)=Y\exp(-\theta'W)$ and $\rho(t)=\log t$.

It is easy to verify that Conditions K1-K2 are satisfied for the above two models if $\Theta\times H$ is assumed to be compact. We will verify Conditions C1-C2 by applying the following Lemma. Let $\psi^{(j)}_\theta(\cdot)$ be $(\partial^j/\partial\theta^j)\psi_\theta(\cdot)$ and $f_{\theta j}(\cdot|z)$ be its conditional density. Denote $f(z)$ as the marginal density of $Z$. Let $M$ be a compact set so that $m_\theta(z)\equiv E[\psi_\theta(Y,W)|Z=z]\in\mbox{int}(M)$ for all $z,\theta$.
\begin{lemma}\label{lem1}
Assume the following conditions hold:
\begin{enumerate}
\item[(a)] $E\{\sup_\theta|\psi_\theta^{(j)}|\}<\infty$ for $j=0,1,2,3$;

\item[(b)] For some even integer $q\geq 10$, $\sup_\theta E\{|\psi_\theta^{(j)}|^q\}<\infty$ for $j=0,1,2,3$;

\item[(c)] $\sup_\theta\sup_{x}|f_{\theta j}^{(r)}(y,w|z)|<\infty$ for $j=0,1,2$ and $r=0,\ldots,4$;

\item[(d)] $\sup_{z}|f^{(r)}(z)|<\infty$ for $r=0,\ldots,4$;

\item[(e)] $0<\inf_z f(z)\leq\sup_z f(z)<\infty$;

\item[(f)] $\sup_{m\in M}|\rho^{(j)}(m)|<\infty$ for $j=0,\ldots,4$.
\end{enumerate}
Suppose that the kernel function $K(\cdot)$ in (\ref{inter03}) satisfies
\begin{eqnarray*}
&&\int K(u)du=1,\;\;\;\int uK(u)du=0,\;\;\;\int u^2K(u)du<\infty,\\
&&\sup_u|K^{(r)}(u)|<\infty\;\;\mbox{for}\;r=0,\ldots,4.
\end{eqnarray*}
Condition C1 holds under the above conditions. If we choose $b_n\asymp n^{-\alpha}$ for
$1/8<\alpha<(q-2)/(4q+16)$, then Condition C2 is satisfied with
\begin{eqnarray}
g&=&2\alpha\wedge\left(\frac{q}{2q+4}-\frac{\alpha(q+4)}{q+2}
-\epsilon\right),\label{relation1}\\
\delta&=&\frac{q}{2q+4}-
\frac{\alpha(2q+6)}{q+2}-2\epsilon\label{relation2}
\end{eqnarray}
for any $\epsilon>0$.
\end{lemma}
The above Lemma specifies the relation between the bandwidth order $\alpha$ in the kernel estimation (\ref{inter03}) and $k^\ast$ in Theorem~\ref{conv-gpl-kmle}. By some algebra, we can verify that $g\in(1/4,1/2]$ and $(2g-1/2)\leq \delta\leq g$ given the above range of $\alpha$ and $q$. We want to point out that the convergence rates of $\widehat\eta(\theta)$ (and its derivatives) may be improved, i.e., larger value of $g$, under more restrictive kernel conditions, see \cite{s89,a95}.

We next apply Theorems~\ref{conv-gpl-kmle}-\ref{gpl-ver} and Lemma~\ref{lem1} to the previous conditional normal (exponential) example, in which $q$ is shown to be arbitrarily large and M is chosen as a sufficiently large compact subset of $(0, \infty)$.
For simplicity, in the below table, we assume that $q=28$, $b_n\asymp n^{-1/5}$, $\epsilon=1/600$ such that $g=151/600>1/4$
and $\delta=1/20$ according to (\ref{relation1})-(\ref{relation2}).
%
\begin{center}
Table 3. {\it Conditional Normal (Exponential) Model $(g=151/600)$}\vspace{0.1in} 
\centering 
\begin{tabular}{c c c}
\hline\hline
 & $\psi=1/2$ & $\psi=1/3$\\
\hline
\mbox{Construction I} & $r_1=1$ & $r_1=2/3$\\
& $k^\ast=1$ & $k^\ast=1$\\
\hline
\mbox{Construction II} & $r_1=451/600$ & $r_1=251/600, r_1=353/600$\\
& $k^\ast=1$ & $k^\ast=2$\\
\multicolumn{3}{@{}p{12.4cm}@{}}{\rule{12.4cm}{0.2pt}}
{\small}
\end{tabular}
\begin{tabular}{c c}
& $\psi=1/4$\\
\hline
\mbox{Construction I} & $r_1=1/2, r_2=1$\\
& $k^\ast=2$\\
\hline
\mbox{Construction II} & $r_1=151/600, r_2=153/600, r_3=157/600, r_4=165/600$\\
 & $r_5=181/600, r_6=213/600, r_7=277/600, r_8=405/600$\\
& $k^\ast=8$ \\
\multicolumn{2}{@{}p{12.6cm}@{}}{\rule{12.6cm}{0.2pt}}\\
\multicolumn{2}{@{}p{12.6cm}@{}}
{\small Remark: $\psi$: convergence rate of $\widehat\theta_n^{(0)}$; $r_k$: Define $\|\widehat\theta_n^{(k)}-\widehat\theta_n\|=O_P(n^{-r_k})$; Construction I: $\widehat I_n$ is constructed by (\ref{semicon1}); Construction II: $\widehat I_n$ is constructed by (\ref{semicon2}).}
\end{tabular}
\end{center}

\subsection{Penalized Estimation in Semiparametric Models}\label{penest}
In many semiparametric models involving a smooth nuisance parameter, it is often convenient and beneficial to perform estimation using
penalization, e.g., \cite{s85, mvg97}. Under regularity conditions, penalized semiparametric log-likelihood estimation can yield fully efficient estimates for $\theta$, see (\ref{penexp}). In penalized estimation framework, the value of $k^\ast$ is shown to relate to the order of the smoothing parameter $\lambda_n$. A surprising result we have is that $k^\ast$ iterations are also sufficient for recovering the estimation sparsity in high dimensional data, see the below partly linear example.

In this subsection, we assume that $\eta$ belongs to the Sobolev class of functions $\mathcal{H}_{k}\equiv\{\eta:J^{2}(\eta)=\int_{\mathcal{Z}}(\eta^{(k)}(z))^{2}dz<\infty\}$, where $\eta^{(j)}$ is the $j$-th derivative of $\eta$ and $\mathcal{Z}$ is some compact set on the real line. The penalized log-likelihood in this context is defined as
\begin{eqnarray}
\log lik_{\lambda_{n}}(\theta,\eta)=n\mathbb{P}_n\log lik(\theta,\eta)-n\lambda_{n}^{2}J^{2}(\eta),\label{penlik}
\end{eqnarray}
where $\lambda_{n}$ is a smoothing parameter. We assume the following bounds for $\lambda_n$:
\begin{eqnarray}
\lambda_{n}=o_{P}(n^{-1/4}) \;\;\mbox{and}\;\;
\lambda_{n}^{-1}=O_{P}(n^{k/(2k+1)}).\label{smooth}
\end{eqnarray}
In practice, $\lambda_{n}$ can be obtained by cross-validation \cite{w98}. Here, the regularized $\widehat S_{n}(\theta)$ becomes the log-profile penalized likelihood $\widehat S_{\lambda_n}(\theta)$:
\begin{eqnarray}
\widehat S_{\lambda_n}(\theta)=\log_{\lambda_n}(\theta,\widehat\eta_{\lambda_n}(\theta)),\label{propen}
\end{eqnarray}
where $\widehat{\eta}_{\lambda_{n}}(\theta)=\arg\sup_{\eta\in\mathcal{H}_{k}}\log lik_{\lambda_{n}}(\theta,\eta)$
for any fixed $\theta$ and $\lambda_{n}$. We define the penalized estimate as $\widehat\theta_{\lambda_n}$.

The construction of the $k$-step penalized estimate $\widehat\theta_{\lambda_n}^{(k)}$ follows from (\ref{parasche-0}) just with the change of $\widehat S_n(\cdot)$ to $\widehat S_{\lambda_n}(\cdot)$. For the penalized estimation, we need to slightly modify Condition G as follows:
\begin{enumerate}
\item[G'.] Assume that, for some constant $c$,
\begin{eqnarray}
\frac{1}{n}\widehat{S}^{(1)}_{\lambda_n}(\theta_0)-c\mathbb{P}_n\widetilde\ell_0&=&O_P(\lambda_n^2),\label{pricon21}\\
\sup_{\theta\in\mathcal{N}(\theta_0)}\left|\frac{1}{n}\widehat{S}^{(2)}_{\lambda_n}
(\theta)+c\widetilde I_0\right|&=&
O_P(\lambda_n\vee\|\theta-\theta_0\|),\label{pricon31}\\
\sup_{\theta\in\mathcal{N}(\theta_0)}\left|\frac{1}{n}\widehat
S_{\lambda_n}^{(3)}(\theta)\right|
&=&O_P(1).\label{pricon41}
\end{eqnarray}
\end{enumerate}
It is easy to verify Condition G' if $\widehat\eta_{\lambda_n}(\theta)$ has an explicit expression and $\log lik_{\lambda_n}(\theta,\eta)$ is smooth w.r.t. $(\theta,\eta)$, see the below example 4. We also want to point out that Condition G' is relaxable to a large extent, see Remark~\ref{remak1}. For example, rather than the explicit form of $\widehat\eta_{\lambda_n}$, we may only require $\widehat\eta_{\lambda_n}$ satisfying  $\|\widehat\eta_{\lambda_n}(\widetilde\theta_n)-\eta_0\|=O_P(\|\widetilde\theta_n-\theta_0\|\vee\lambda_n)$ for any consistent $\widetilde\theta_n$.

In view of (\ref{lfs}) and (\ref{gpl-info2}), we can prove Theorem~\ref{cor-pen} similarly as Theorem~\ref{conv-gpl-kmle}.
\begin{theorem}\label{cor-pen}
Suppose Condition G' holds, the penalized MLE $\widehat\theta_{\lambda_n}$ is consistent and $\widetilde I_0$ is nonsingular. We have
\begin{eqnarray}
\sqrt{n}(\widehat{\theta}_{\lambda_n}-\theta_{0})=
\frac{1}{\sqrt{n}}\sum_{i=1}^{n}\widetilde{I}_{0}^{-1}\widetilde{\ell}_{0}
(X_{i})+O_P(\sqrt{n}\lambda_n^2).\label{penexp}
\end{eqnarray}
Define $g=\max\{g': \lambda_n=O_P(n^{-g'})\}$, and thus $1/4<g\leq k/(2k+1)$ based on Condition (\ref{smooth}). Construct $\widehat\theta_{\lambda_n}^{(k)}$ as in (\ref{parasche-0}) with the change of $\widehat S_n(\cdot)$ to $\widehat S_{\lambda_n}(\cdot)$. Then all the conclusions for $\widehat\theta_n^{(k)}$ in Theorem~\ref{conv-gpl-kmle} also hold for $\widehat\theta_{\lambda_n}^{(k)}$.
\end{theorem}
The above asymptotic linear expansion (\ref{penexp}) was also derived in \cite{ck09} but under very different conditions. Theorem~\ref{cor-pen} implies that $k^\ast$ depends on the order of the smoothing parameter $\lambda_n$, i.e., the value of $g$, see (\ref{gpl-thm-2}). Because of the duality between the penalized estimation and sieve estimation, we expect that the above conclusions also hold for the semiparametric sieve estimation, see \cite{c07}. For example, when $\eta_0$ is estimated in the form of B-spline (local polynomial) as in \cite{hzz07} (\cite{fhw95,cfgw97}) , $k^\ast$ may rely on the growth rate of the number of basis functions (the order of bandwidth in the kernel function). The detailed theoretical exploration towards this direction is not considered in this article due to the space limitation.

We next apply Theorem~\ref{cor-pen} to the following partly linear models under high dimensional data. Surprisingly, we discover that one step iteration is sufficient for achieving the semiparametric estimation efficiency and recovering the estimation sparsity simultaneously.

{\it Example 4. Sparse and Efficient Estimation of Partial Spline Model}

The partial smoothing spline represents an important class of semiparametric models under penalized estimation. In particular, we consider
\begin{eqnarray}
Y=W'\theta+\eta(Z)+\epsilon,\label{partly}
\end{eqnarray}
where $\eta\in\mathcal{H}_k$ and $0\leq Z\leq 1$. For simplicity, we assume that $\epsilon\overset{iid}{\sim}N(0,\sigma^2)$ and is independent of $(W,Z)$. The normality of $\epsilon$ can be relaxed to the sub-exponential tail condition. In this example, we assume that some components of $\theta_0$ are exactly zero which is common for high dimensional data. It is well known that effective variable selection in semiparametric models could greatly improve their prediction accuracy and interpretability, e.g., \cite{bm05,cz10}. To achieve the estimation efficiency and recover sparsity of $\theta$, Cheng and Zhang (2010) proposed the following double penalty estimation approach for (\ref{partly}). Specifically, they define $(\widehat\theta_{\lambda_n},\widehat\eta_{\lambda_n})$ as the minimizer of
\begin{eqnarray}
n\mathbb{P}_n(Y-W'\theta-\eta(Z))^2+n\lambda_n^2 J^2(\eta)+n\tau_n^2\sum_{j=1}^{d}\frac{|\theta_j|}{|\widetilde\theta_j|^\gamma},\label{psp}
\end{eqnarray}
where $\gamma$ is a fixed positive constant and $\widetilde\theta=(\widetilde\theta_1,\ldots,\widetilde\theta_d)'$ is the consistent initial estimate, over $\Theta\times\mathcal{H}_k$. 

We will show that $\widehat\theta_{\lambda_n}^{(1)}$ possesses the same {\it semiparametric oracle property}, whose definition is given below, as $\widehat\theta_{\lambda_n}$. The standard smoothing spline theory suggests that
\begin{eqnarray}
\widehat\eta_{\lambda_n}(\theta)({\zv})=A(\lambda_n)(\yv-\wv\theta),\label{pssnon}
\end{eqnarray}
where $\widehat\eta_{\lambda_n}(\theta)(\zv)=(\widehat\eta_{\lambda_n}(\theta)(z_1),\ldots,\widehat\eta_{\lambda_n}(\theta)(z_n))'$, $\yv=(y_1,\ldots,y_n)'$ and $\wv=(w_1',\ldots,w_n')'$. The expression of the $n\times n$ influence matrix $A(\lambda_n)$ can be found in \cite{h86}. Therefore, $\widehat\eta_{\lambda_n}(\theta)$ is a natural spline of order $(2k-1)$ with knots on $z_i$'s for any fixed $\theta$. Plugging (\ref{pssnon}) back to (\ref{psp}), we have
\begin{eqnarray}
\widehat S_{\lambda_n}(\theta)=\widetilde S_{\lambda_n}(\theta)+n\tau_n^2\sum_{j=1}^{d}\frac{|\theta_j|}{|\widetilde\theta_j|^\gamma},\label{slambdan}
\end{eqnarray}
where
\begin{eqnarray}
\widetilde S_{\lambda_n}(\theta)=(\yv-\wv\theta)'[I-A(\lambda_n)](\yv-\wv\theta)\label{tildeslam}
\end{eqnarray}
and $I$ is the identity matrix of size $n$. When $\tau_n=0$,
the minimizer of (\ref{psp}) becomes the partial smoothing spline, and we denote it as $(\widetilde\theta_{\lambda_n},\widetilde\eta_{\lambda_n})$.
Note that $\widetilde\theta_{\lambda_n}$ has a simple analytic form as $\widetilde\theta_{\lambda_n}=[\wv'(I-A(\lambda_n))\wv]^{-1}\wv'[I-A(\lambda_n)]\yv$. However, $\widehat\theta_{\lambda_n}$ as the minimizer of $\widehat S_{\lambda_n}(\theta)$ does not have an explicit solution form, and has to be iteratively computed using software like Quadratic Programming or LARS \cite{ehjt04}, see Section 4 of \cite{cz10}. Specifically, based on (\ref{parasche-0})-(\ref{semicon1}), we construct $\widehat\theta_{\lambda_n}^{(1)}$ as follows:
\begin{eqnarray*}
\widehat\theta_{\lambda_n}^{(1)}=\widehat\theta_{\lambda_n}^{(0)}+\left[\frac{\wv'(I-A(\lambda_n))\wv}{n}\right]^{-1}\left[
\frac{\wv'(I-A(\lambda_n))(\yv-\wv\widehat\theta_{\lambda_n}^{(0)})}{n}-\frac{\tau_n^2}{2}\delta_n(\widehat\theta_{\lambda_n}^{(0)})\right],
\end{eqnarray*}
where $\delta_n(\theta)=(sign(\theta_1)/|\widetilde\theta_1|^\gamma,\ldots,sign(\theta_d)/|\widetilde\theta_d|^\gamma)'$.

Without loss of generality, we write $\theta_0=(\theta_1',\theta_2')'$, where $\theta_1$ consists of all $q$ nonzero components and $\theta_2$ consists of the rest $(d-q)$ zero elements, and define $\widehat\theta_{\lambda_n}=(\widehat\theta_{\lambda_n,1}'\widehat\theta_{\lambda_n,2}')'$ accordingly. We assume that $W$ has zero mean, strictly positive definite covariance matrix $\Sigma$ and finite fourth moment. The observations $z_i$'s (real numbers) are sorted and satisfy
\begin{eqnarray}
\int_{0}^{z_i}u(w)dw=\frac{i}{n}\;\;\;\;\;\;\;\;\;\mbox{for}\;i=1,2,\ldots,n,\label{recon1}
\end{eqnarray}
where $u(\cdot)$ is a continuous and strictly positive function. The above regularity conditions are commonly used in the literature, e.g., \cite{h86,cw79}, and are relaxable. For example, Condition (\ref{recon1}) can be weakened to the case in which $z_i$'s are sufficiently close to a sequence satisfying (\ref{recon1}). For simplicity, we assume that $\gamma=1$ and $\widetilde\theta$ is $\sqrt{n}$-consistent. In this example, $\widetilde\theta_{\lambda_n}$ or the difference based estimate \cite{y97}, which are both known to be $\sqrt{n}$ consistent, can serve as $\widetilde\theta$ or $\widehat\theta_{\lambda_n}^{(0)}$.

In this example, we say $\widehat\theta_{\lambda_n}$ satisfies the {\it semiparametric oracle property} if
\begin{enumerate}
\item[O1.] $\sqrt{n}(\widehat{\theta}_{\lambda_n,1}-\theta_{1})\overset{d}{\longrightarrow}
N(0,\sigma^2\Sigma_{11}^{-1})$, where $\Sigma_{11}$ is the $q\times q$ upper-left submatrix of $\Sigma$ [Semiparametric Efficiency];

\item[O2.] $\widehat\theta_{\lambda_n,2}=0$ with probability tending to one [Sparsity].
\end{enumerate}
It is easily shown that $\sigma^2\Sigma_{11}^{-1}$ in O1 is the semiparametric efficiency bound for $\theta_1$ since $z$ is assumed to be fixed.
\begin{corollary}\label{sparseco}
If $n^{k/(2k+1)}\lambda_n\rightarrow\lambda_0>0$ and $n^{k/(2k+1)}\tau_n\rightarrow\tau_0>0$, then $\widehat\theta_{\lambda_{n}}$ is $\sqrt{n}$-consistent and satisfies the semiparametric oracle property. Given that $\widehat\theta_{\lambda_n}^{(0)}$ is $\sqrt{n}$-consistent, then $\|\widehat\theta_{\lambda_n}^{(1)}-\widehat\theta_{\lambda_n}\|=O_P(n^{-1})$ and $\widehat\theta_{\lambda_n}^{(1)}$ also enjoys the semiparametric oracle property.
\end{corollary}
The above Corollary is a simple but interesting application of Theorem~
\ref{cor-pen}. We can definitely relax its conditions to the general $\gamma$ and non-$\sqrt{n}$ consistent $\widehat\theta_{\lambda_n}^{(0)}$ in which we may require more than one iteration. The conditions on $\lambda_n$ and $\tau_n$ are also chosen for simplicity of expositions and are relaxable. In addition, the proof of Corollary~\ref{sparseco} implies the following special case of (\ref{penexp}):
\begin{eqnarray*}
\sqrt{n}(\widehat{\theta}_{\lambda_n,1}-\theta_{1})=\frac{1}{\sqrt{n}}\Sigma_{11}^{-1}\sum_{i=1}^{n}W_{1i}\epsilon_i
+O_P\left(n^{-\frac{2k-1}{2(2k+1)}}\right),
\end{eqnarray*}
where $W_{1i}$ is the first $q$ elements of $W_i$. It is also possible to extend the conclusions of Corollary~\ref{sparseco} to the semiparametric quasi-likelihood framework proposed in \cite{mvg97} after more tedious algebra.

\section{Initial Estimate}\label{iniest}
In this paper, we assume the existence of a $n^\psi$-consistent $\widehat{\theta}_n^{(0)}$ just as the numerical result assumes the iterations commence in a suitable neighborhood of $\theta_0$. Occasionally, the semiparametric model structure can be exploited to produce a $\sqrt{n}$-consistent initial estimate, e.g., \citep{y97, s96}. However, if such ad-hoc methods are unavailable, a general strategy is to conduct a search of some objective function at finitely many $\theta$-value and call the optimizer as the initial estimate. The numerical analysis literature suggest several search
strategies for parametric models, e.g. \citep{s72,f80}, and Robinson (1988) subsequently proved the consistency and convergence rate of those numerical outcome. In this section, we extend Robinson's results to semiparametric models, i.e., Theorem~\ref{init}. This extension is nontrivial since our objective function usually has no explicit form and is possibly nonsmooth. In fact, our theoretical results on searching $\widehat\theta_n^{(0)}$ can be applied to any objective functions satisfying the below Conditions I1-I2, and are thus of independent interest.

We use the generalized profile likelihood $\widehat S_n(\theta)$ as our objective function in semiparametric models. Besides the compactness of $\Theta$ and consistency of $\widehat\theta_n$, we have two primary conditions I1-I2 on $\widehat S_n(\theta)$ to guarantee the validity of the grid search methods we will consider.
\begin{itemize}
\item[I1.] [Asymptotic Uniqueness] For any random sequence $\{\widetilde{\theta}_{n}\}\in\Theta$,
\begin{eqnarray}
\;\;\;\;\;\;\;[\widehat S_n(\widetilde\theta_n)-\widehat S_n(\widehat\theta_n)]/n=o_{P}(1)\;\;\mbox{implies that}\;\;\widetilde{\theta}_{n}-\theta_{0}
=o_{P}(1).\label{asyuni}
\end{eqnarray}

\item[I2.] [Asymptotic Expansion] For any consistent $\widetilde\theta_n$, $\widehat S_n$ satisfies
\begin{eqnarray}
\widehat S_n(\widetilde\theta_n)&=&\widehat S_n(\theta_0)+n(\widetilde\theta_n-\theta_0)'\mathbb{P}_n\widetilde\ell_0-\frac{n}{2}
(\widetilde\theta_n-\theta_0)'\widetilde I_0(\widetilde\theta_n-\theta_0)\nonumber\\&&+\Delta_{n}(\widetilde\theta_{n}),\label{iniass1}
\end{eqnarray}
where $\Delta_{n}(\theta)=n\|\theta-\theta_0\|^3\vee n^{1-2r}\|\theta-\theta_0\|$ and $1/4<r\leq 1/2$.
\end{itemize}
Condition I1 is usually implied by the model identifiability conditions. Note that, in Condition I2, we only assume the existence of the asymptotic expansion (\ref{iniass1}) but not assume the continuity of $\widehat S_n(\cdot)$. In Section~\ref{semisec}, we have shown that the log-profile likelihood $\log pl_n(\cdot)$ as a special case of $\widehat S_n(\cdot)$ satisfies I2 under model Assumptions M1-M4, see (\ref{lnplexp}). As for the regularized $\widehat S_n(\cdot)$, we can verify I2 under Condition G using a three term Taylor expansion of $\widehat S_n$. Specifically, I2 is satisfied if we assume Conditions K1-K2 \& C1-C2 (G') for the kernel estimation (penalized estimation). In particular, we can change $n^{1-2r}$ to $n\lambda_n^2$ in $\Delta_{n}(\cdot)$ when considering the penalized estimation.

Now we consider two types of grid search: deterministic type and stochastic type. In the former, we form a grid of cubes with sides of length $sn^{-\psi}$ over $\mathbb{R}^d$ for some $s>0$ and $0<\psi\leq 1/2$, and thus obtain a set of points $\mathcal{D}_n=\{\theta_{iD}\}$ regularly spaced throughout $\Theta$ with cardinality $card(\mathcal{D}_n)\geq Cn^{d\psi}$ for some $C>0$. The grid point which maximizes $\widehat S_n(\theta)$ is thought of as $\widehat{\theta}_{n}^{(0)}$. However, this deterministic search could be very slow if the dimension $d$ of $\theta$ is high. This motivates us to propose the stochastic search in which the search points are the realizations of some independent random variable $\bar\theta$ with strictly positive density around $\theta_0$, e.g., $\bar\theta\sim Unif[\Theta]$. And we require that the magnitude of the stochastic search points remains $n^{\psi}$ no matter how large the dimension $d$ is. In theory, the stochastic grid search has significant computational savings over the deterministic alternative. In the below Theorem~\ref{init} we rigorously prove that the convergence rates of the above numerical outcomes are $n^{\psi}$-consistent for $0<\psi\leq1/2$.

\begin{theorem}\label{init}
Let $\mathcal{D}_{n}$ be a set of points regularly spaced throughout $\Theta$ with $card(\mathcal{D}_n)\geq Cn^{d\psi}$ for some $C>0$ and $0<\psi\leq1/2$. Assume that $\bar{\theta}$ is independent of the data and admits a density having support $\Theta$ and bounded away from zero in some neighborhood of $\theta_{0}$. Let $\mathcal{S}_{n}$ be a set of realizations of $\bar{\theta}$ with $card(\mathcal{S}_n)\geq\widetilde Cn^{\psi}$ for some $\widetilde C>0$ and $0<\psi\leq 1/2$.
Suppose that Conditions I1-I2 hold, and that the parameter space $\Theta$ is compact. Then, if $\widehat\theta_n$ defined in (\ref{semiest}) is consistent and $\widetilde I_0$ is nonsingular, we have
\begin{eqnarray}
\theta_{n}^{D}-\theta_{0}&=&O_{P}(n^{-\psi}),\label{thm5}\\
\theta_{n}^{S}-\theta_{0}&=&O_{P}(n^{-\psi}),\label{thm6}
\end{eqnarray}
where $\theta_{n}^{D}=\arg\max_{\theta\in\mathcal{D}_{n}}\widehat S_{n}(\theta)$ and $\theta_{n}^{S}=\arg\max_{\theta\in\mathcal{S}_{n}}\widehat S_{n}(\theta)$.
\end{theorem}
Theorems~\ref{thm1}-\ref{conv-gpl-kmle} together with the above Theorem~\ref{init} offer rigorous statistical analysis for the
general iterative semiparametric estimation algorithm presented in Introduction section. Those theorems indicate a tradeoff
between the computational cost of searching for an initial estimate, i.e. $card(\mathcal{D}_n)$ or $card(\mathcal{S}_n)$, and that of generating
an efficient estimate, i.e., $k^\ast$. Theorem~\ref{init} can be applied to all the examples we have considered. Specifically,
Condition I1 is verified for Examples 1-2 in \cite{ck08b}, and we can easily verify Condition I1 in Example 3 by adapting the consistency
proof of $\widehat\theta_n$ in  \cite{sw92}, see its Proposition 1. In fact, the Conditions I1-I2 are very mild and can be satisfied in a wide
range of semiparametric models, e.g., proportional odds model and penalized semiparametric logistic regression.
%

\vskip 1em \centerline{\Large \bf APPENDIX} \vskip -2em
\setcounter{subsection}{0}
\renewcommand{\thesubsection}{A.\arabic{subsection}}
\setcounter{equation}{0}
\renewcommand{\theequation}{A.\arabic{equation}}
\setcounter{lemma}{0}
\renewcommand{\thelemma}{A.\arabic{lemma}}\hspace{0.5in}

\subsection{Conditions M1-M4 on the Least Favorable Submodel}
The LFS in Section~\ref{semisec} is constructed in
the following manner. We first assume the existence of a smooth
map from the neighborhood of $\theta$ into $\mathcal{H}$, of the form $t\mapsto\eta_\ast(t;\theta, \eta)$, such that
the map $t\mapsto \ell(t,\theta,\eta)(x)$ can be defined as follows:
\begin{eqnarray}
\ell(t,\theta,\eta)(x)&=&\log lik(t,\eta_{\ast}(t;\theta,\eta))(x),
\end{eqnarray}
where we require $\eta_{\ast}(\theta;\theta,\eta)=\eta$ for
all $(\theta,\eta)\in\Theta\times{\cal H}$. Thus, $\log pl_n(\theta)=\sum_{i=1}^{n}\ell(X_i;\theta,\theta,\widehat{\eta}(\theta))$.
See \cite{ck08b} for similar constructions. We define $\dot{\ell}(t,\theta,\eta)$,
$\ddot{\ell}(t,\theta,\eta)$ and $\ell^{(3)}(t,\theta,\eta)$ as the first, second and third derivative of
$\ell(t,\theta,\eta)$ with respect to $t$, respectively. Also denote
$\ell_{t,\theta}(t,\theta,\eta)$ as $(\partial^2/\partial t\partial\theta)\ell(t,\theta,\eta)$.
\begin{enumerate}
\item[M1.]We assume that the derivatives $(\partial^{l+m}/\partial t^l\partial\theta^m)\ell(t,\theta,\eta)$ have
integrable envelop functions in $L_1(P)$ for $(l+m)\leq 3$, and that the Fr\'{e}chet
derivatives of $\eta\mapsto\ddot\ell(\theta_0,\theta_0,\eta)$ and $\eta\mapsto\ell_{t,\theta}(\theta_0,\theta_0,\eta)$ are bounded around
$\eta_0$;

\item[M2.] $E\dot\ell(\theta_0,\theta_0,\eta)=O(\|\eta-\eta_0\|^2)$ for all $\eta$ around $\eta_0$;

\item[M3.] $\mathbb{G}_n(\dot{\ell}(\theta_{0},\theta_0,
\widehat\eta(\widetilde\theta_n))-\dot{ \ell}(\theta_0,\theta_0,\eta_0))=O_P(n^{-2r+1/2}\vee n^{1/2-r}\|\widetilde\theta_n-\theta_0\|)$
for any $\widetilde\theta_n\overset{P}{\rightarrow}\theta_0$;

\item[M4.] The classes of functions
$\{\ddot{\ell}(t,\theta,\eta)(x): (t,\theta,\eta)\in V\}$ and
$\{\ell_{t,\theta}(t,\theta,\eta)(x): (t,\theta, \eta)\in V\}$ are
$P$-Donsker, and $\{\ell^{(3)}(t,\theta,\eta)(x):
(t,\theta,\eta)\in V\}$ is $P$-Glivenko-Cantelli, where $V$ is some neighborhood of $(\theta_{0}, \theta_{0},
\eta_{0})$.
\end{enumerate}
See Section 2.2 of \citep{ck08b} for the discussions on M1-M4.

\subsection{Useful Lemmas}
The first two Lemmas are used in the proof of Lemma~\ref{step}. The Lemmas~\ref{gpl-le3}, \ref{lemm-para}, \ref{gpl-le1} and \ref{sparlem} are used in the proofs of Theorem~\ref{conv-gpl-kmle}, Corollary~\ref{conv-para}, Theorem~\ref{gpl-ver} and Corollary~\ref{sparseco}, respectively.
\begin{lemma}\label{mainthm}
Suppose that Conditions M1-M4 and (\ref{convrate}) hold. If $\widetilde{\theta}_n$ is $n^{\psi}$-consistent, then we have
\begin{eqnarray}
\;\;\;\;\;\;\;\widehat{\ell}_{n}(\widetilde{\theta}_n,s_{n})&=&\mathbb{P}_{n}\widetilde{\ell}_{0}+O_{P}\left(n^{-\psi}\vee
|s_{n}|\vee\frac{g_{r}(n^{-\psi}\vee|s_{n}|)}
{n|s_{n}|}\right),\label{ga1}\\
\widehat{\ell}_{n}(\widehat{\theta}_n+U_{n},s_{n})&=&\widehat{\ell}_{n}
(\widehat{\theta}_n,s_{n})-\widetilde{I}_{0}U_{n}\nonumber\\&&+O_{P}\left(\frac{g_r(|s_n|\vee \|U_n\|)\vee n^{1/2-2r}}{n|s_n|}\right),\label{lem2}\\
\widehat{I}_{n}(\widetilde{\theta}_{n},t_{n})&=&\widetilde{I}_{0}\nonumber\\&&+O_{P}
\left(\frac{g_{r}(\|\widetilde\theta_n-\widehat\theta_n\|\vee
|t_{n}|)\vee nt_n\|\widetilde\theta_n-\widehat\theta_n\|\vee n^{1/2-2r}}{nt_{n}^{2}}\right),\label{estrel}
\end{eqnarray}
where $g_r(t)=nt^3\vee
n^{1-2r}t$ and $U_{n}=O_{P}(n^{-s})$ for some $s>0$.
\end{lemma}

{\sc Proof: }
Under the assumptions M1-M4 and (\ref{convrate}), \citep{ck08b} proved the following asymptotic expansion of $\log pl_n(\bar\theta_n)$, where $\bar\theta_n$ is consistent,
\begin{eqnarray}
\log pl_{n}(\bar{\theta}_{n})&=&\log pl_{n}(\theta_{0})+(\bar{\theta}_{n}-\theta_{0})'
\sum_{i=1}^{n}\widetilde{\ell}_{0}(X_{i})-\frac{n}{2}(\bar{\theta}_
{n}-\theta_{0})'\widetilde{I}_{0}(\bar{\theta}_{n}-\theta_{0})
\nonumber\\&&+O_P\left(g_{r}(\|\bar{\theta}_n-\theta_0\|)\right),\label{lnplexp}\\
\log pl_{n}(\bar{\theta}_{n})&=&\log
pl_{n}(\widehat{\theta}_{n})-\frac{1}{2}n(\bar{\theta}_{n}-
\widehat{\theta}_{n})'\widetilde{I}_{0}(\bar{\theta}_{n}-
\widehat{\theta}_{n})\nonumber\\&&+O_{P}\left(g_{r}(\|\bar\theta_n-\widehat\theta_n\|)\vee
n^{1/2-2r}\right).\label{lnplexphat}
\end{eqnarray}
We first prove (\ref{lem2}).
(\ref{lnplexphat}) implies
that
\begin{eqnarray*}
\log pl_{n}(\widehat\theta_n+V_{n}+s_{n}v_{i})&=&\log
pl_{n}(\widehat\theta_{n})-\frac{n}{2}(V_{n}+s_{n}v_{i})'\widetilde{I}_{0}(V_{n}+s_{n}v_{i})\\
&&+O_P(g_{r}(|s_{n}|\vee\|V_{n}\|)\vee n^{1/2-2r}),\\
\log pl_{n}(\widehat\theta_n+V_{n})&=&\log
pl_{n}(\widehat\theta_{n})-\frac{n}{2}V_{n}'\widetilde{I}_{0}V_{n}+O_P(g_{r}(\|V_{n}\|)\vee n^{1/2-2r}),
\end{eqnarray*}
for any random vector $V_{n}=o_{P}(1)$ and
$s_{n}\overset{P}{\rightarrow}0$. Combining the above two expansions and (\ref{estesco}), we have
$$[\widehat\ell_{n}(\widehat\theta_n+V_n,s_n)]_i=-\frac{s_n}{2}v_i'\widetilde I_0v_i-v_i'\widetilde I_0 V_n+O_P\left(\frac{g_r(|s_n|\vee\|V_n\|)\vee n^{1/2-2r}}{n|s_n|}\right).$$
By taking $V_n=0$ and $U_n$, respectively, in the above equation, we have proved (\ref{lem2}). Following similar analysis in the above, (\ref{estesco}) \& (\ref{lnplexp}) yield (\ref{ga1}), and (\ref{estei}) \& (\ref{lnplexphat}) yield (\ref{estrel}). This completes the whole
proof. $\Box$
\begin{lemma}\label{mainthm2}
Suppose that Conditions M1-M4 and (\ref{convrate}) hold. If
\begin{eqnarray}
\widehat I_{n}(\widehat{\theta}_{n}^{(k-1)},t_{n})-\widetilde{I}_{0}=O_{P}(r_{n}^{(k-1)}),\label{ass-pi}
\end{eqnarray}
then we have $\|\widehat{\theta}_{n}^{(k)}-\widehat{\theta}_{n}\|=$
\begin{eqnarray}
&&O_{P}\left(|s_{n}|\vee\|\widehat{\theta}_{n}^{(k-1)}-\widehat{\theta}_{n}\|r_{n}^{(k-1)}\vee\right.\nonumber\\
&&\mbox{\hspace{1.5in}}\left.
\frac{g_{r}(|s_{n}|\vee\|\widehat{\theta}_{n}^{(k-1)}-\widehat{\theta}_{n}\|)\vee n^{1/2-2r}}{n|s_{n}|} \right)\label{res1}
\end{eqnarray}
for $k=1,2,\ldots$.
\end{lemma}
{\sc Proof:} Based on (\ref{sche}), we have
\begin{eqnarray}
\widehat I_{n}(\widehat{\theta}_{n}^{(k-1)},t_{n})\sqrt{n}(\widehat{\theta}_{n}^{(k)}-\widehat{\theta}_{n})&=&
\left[\sqrt{n}\widehat I_{n}(\widehat{\theta}_{n}^{(k-1)},t_{n})(\widehat{\theta}_{n}^{(k-1)}-\widehat{\theta}_{n})\right]+\sqrt{n}\widehat\ell_{n}
(\widehat{\theta}_{n},s_{n})\nonumber\\&&+
\left[\sqrt{n}(\widehat\ell_{n}(\widehat{\theta}_{n}^{(k-1)},s_{n})-\widehat\ell_{n}
(\widehat{\theta}_{n},s_{n}))\right].\label{pro1}
\end{eqnarray}
The second term in (\ref{pro1}) equals to
\begin{eqnarray*}
O_{P}\left(\sqrt{n}|s_{n}|\vee \frac{g_{r}(|s_{n}|)\vee
n^{1/2-2r}}{\sqrt{n}|s_{n}|}\right)
\end{eqnarray*}
according to (\ref{estesco}) and (\ref{lnplexphat}). The third
term in (\ref{pro1}) can be written as
\begin{eqnarray*}
-\sqrt{n}\widetilde{I}_{0}(\widehat{\theta}_{n}^{(k-1)}-\widehat{\theta}_{n})+O_{P}\left(\frac{g_{r}
(|s_{n}|\vee\|\widehat{\theta}_{n}^{(k-1)}-\widehat{\theta}_{n}\|)\vee n^{1/2-2r}}{\sqrt{n}|s_{n}|}\right).
\end{eqnarray*}
for $k=1,2,\ldots$ by replacing $U_{n}$ with
$(\widehat{\theta}_{n}^{(k-1)}-\widehat{\theta}_{n})$ in
(\ref{lem2}). Combining the above analysis, the assumption
(\ref{ass-pi}) and nonsingularity of $\widetilde{I}_{0}$, we complete
the proof of (\ref{res1}). $\Box$

\begin{lemma}\label{gpl-le3}
Suppose that Condition G holds. If $\widetilde{\theta}_n$ is a $n^{\psi}$-consistent
estimator for $0<\psi\leq 1/2$, then we have
\begin{eqnarray}
&&n^{-1}[\widehat{S}_n^{(1)}(\widetilde{\theta}_n)-
\widehat{S}_n^{(1)}(\theta_0)]\nonumber\\&=&-
\widetilde{I}_0(\widetilde{\theta}_n-\theta_0)+O_P((n^{-g}\vee\|\widetilde{\theta}_n-\theta_0
\|)\|\widetilde{\theta}_n-\theta_0
\|),\label{gpl-kel-eq0}\\
&&n^{-1}[\widehat{S}_n^{(1)}(\widetilde{\theta}_n+U_n)-
\widehat{S}_n^{(1)}(\widetilde{\theta}_n)]\nonumber\\&=&-
\widetilde{I}_0U_n+O_P((n^{-g}\vee\|\widetilde{\theta}_n
-\theta_0\|)\|U_n\|),\label{gpl-kel-eq1}
\end{eqnarray}
where $U_n$ a statistic of the order $O_P(n^{-s})$ for some $s\geq
\psi$.
\end{lemma}

{\sc Proof:}
We first consider (\ref{gpl-kel-eq0}). Using a Taylor's expansion,
we have
\begin{eqnarray*}
\frac{1}{n}\widehat{S}_n^{(1)}(\widetilde{\theta}_n)&=&
\frac{1}{n}\widehat{S}_n^{(1)}(\theta_0)+\frac{1}{n}
\widehat{S}_n^{(2)}(\theta_0)(\widetilde{\theta}_n-\theta_0)\nonumber\\&&+
\frac{1}{2}(\widetilde{\theta}_n-\theta_0)\otimes\frac{
\widehat{S}_n^{(3)}(\theta^\ast_1)}{n}\otimes
(\widetilde{\theta}_n-\theta_0)\nonumber\\
&=&\frac{1}{n}\widehat{S}_n^{(1)}(\theta_0)+A+B,
\end{eqnarray*}
where $\theta^\ast_1$ lies between $\widetilde{\theta}_n$ and
$\theta_0$. In view of (\ref{gpl-info2}) and
(\ref{pricon3}), we have
$A=-\widetilde{I}_0(\widetilde{\theta}_n-\theta_0)+O_P(
n^{-g}\|\widetilde{\theta}_n-\theta_0\|)$. Condition (\ref{pricon4}) implies that
$B=O_P(\|\widetilde{\theta}_n-\theta_0\|^2)$. This completes the
proof of (\ref{gpl-kel-eq0}).
We next consider (\ref{gpl-kel-eq1}). Similarly, we have $[\widehat
S_n^{(1)}(\widetilde{\theta}_n+U_n)-\widehat
S_n^{(1)}(\widetilde{\theta}_n)]/n$
\begin{eqnarray*}
&=&\frac{1}{n}\widehat
S_n^{(2)}(\widetilde{\theta}_n)U_n+O_P(\|U_n\|^2)\\
&=&\frac{1}{n}S_n^{(2)}
(\widetilde{\theta}_n)U_n+O_P(n^{-g}\|U_n\|\vee\|U_n\|^2),\\
&=&\frac{1}{n}S_n^{(2)}(\theta_0)U_n+O_P(\|\widetilde{\theta}_n-\theta_0\|\|U_n\|
\vee n^{-g}\|U_n\|\vee\|U_n\|^2),\\
&=&-\widetilde{I}_0U_n+O_P(n^{-1/2}\|U_n\|\vee\|\widetilde{\theta}_n-\theta_0\|\|U_n\|\vee
n^{-g}\|U_n\|\vee\|U_n\|^2),
\end{eqnarray*}
where the second equation follows from (\ref{pricon3}), the third
equality follows from (\ref{pricon4}) and the last equation follows
from CLT and (\ref{gpl-info2}). Considering that $1/4<g\leq 1/2$
and $s\geq \psi$, we have proved (\ref{gpl-kel-eq1}). $\Box$

\begin{lemma}\label{lemm-para}
Let $\widehat S_n(\theta)=\sum_{i=1}^{n}\ell_\theta(X_i)$. Suppose that $\widetilde{\theta}_n$ is a $n^\psi$-consistent
estimator for $0<\psi\leq 1/2$. If $\ell_\theta(\cdot)$ satisfies
P1 \& P2, we have
\begin{eqnarray}
&&n^{-1}[\widehat{S}_n^{(1)}(\widetilde{\theta}_n)-
\widehat{S}_n^{(1)}(\theta_0)]\nonumber\\&=&-
I_0(\widetilde{\theta}_n-\theta_0)+O_P(\|\widetilde{\theta}_n-\theta_0
\|^2),\label{gpl-kel-eq0-para}\\
&&n^{-1}[\widehat{S}_n^{(1)}(\widetilde{\theta}_n+U_n)-
\widehat{S}_n^{(1)}(\widetilde{\theta}_n)]\nonumber\\&=&-
I_0U_n+O_P(\|\widetilde{\theta}_n
-\theta_0\|\|U_n\|),\label{gpl-kel-eq1-para}
\end{eqnarray}
where $U_n$ a statistic of the order $O_P(n^{-s})$ for any $s\geq
\psi$.
\end{lemma}
{\sc Proof:} We only provide the proof of (\ref{gpl-kel-eq1-para}) since that of
(\ref{gpl-kel-eq0-para}) is completely analogous and simpler. To show
(\ref{gpl-kel-eq1-para}), it suffices to prove that, for every $C_1, C_2>0$
and $s\geq\psi$,
\begin{eqnarray*}
&&\sup_{|t|\leq C_1, |u|\leq
C_2}\left|n^{-1}\left[\widehat S_n^{(1)}(\theta_0+n^{-\psi}t+n^{-s}u)
-\widehat S_n^{(1)}(\theta_0+n^{-\psi}t)\right]+n^{-s}I_0u\right|
\\&=&O_P(n^{-s-\psi}).
\end{eqnarray*}
Denote $Z_n(t,u)=n^{-1/2}[\widehat S_n^{(1)}(\theta_0+n^{-\psi}t+n^{-s}u)-
\widehat S_n^{(1)}(\theta_0+n^{-\psi}t)]$ and $Z_n^0(t,u)=Z_n(t,u)-EZ_n(t,u)$.
Then, it suffices to show that
\begin{eqnarray}
\sup_{|t|\leq C_1, |u|\leq C_2}|Z_n^0(t,u)|&=&O_P(n^{1/2-\psi-s}),
\label{inter4}\\
\sup_{|t|\leq C_1, |u|\leq
C_2}|EZ_n(t,u)+n^{1/2-s}I_0u|&=&O_P(n^{1/2-\psi-s}).\label{inter50}
\end{eqnarray}
The proofs of (\ref{inter4}) and (\ref{inter50}) are similar as
those of (2.3) and (2.4) in Page 1224 of \citep{jjv85}, and are
thus skipped. $\Box$

\begin{lemma}\label{gpl-le1}
Suppose Conditions K1-K2 \& C1-C2 hold. Then we have
\begin{eqnarray}
\frac{1}{\sqrt{n}}\sum_{i=1}^{n}\left(\frac{\partial}{\partial\theta}|_{\theta=\theta_0}
A_{\theta,\eta_\ast(\theta)}\right)[\widehat{\eta}_0-\eta_{\ast0}]
(X_i)&=&O_P(n^{-\delta}),\label{gpl-inter1}\\
\frac{1}{\sqrt{n}}\sum_{i=1}^{n}A_{\theta_0,\eta_0}[\widehat{\eta}_0^{(1)}
-\eta_{\ast0}^{(1)}](X_i)&=&O_P(n^{-\delta}),\label{gpl-inter2}\\
\frac{1}{\sqrt{n}}\dot{r}_{n}(\theta_0)&=&O_P(n^{1/2-2g}),\label{lem-ine3}
\end{eqnarray}
where
$r_{n}(\theta)\equiv\widehat{S}_n(\theta)-S_n(\theta)-\sum_{i=1}^{n}A_{\theta,\eta_\ast(\theta)}
[\widehat{\eta}(\theta)-\eta_\ast(\theta)]$.
\end{lemma}
{\sc Proof:} The proof of Lemma 2 in
\citep{sw92} directly implies (\ref{gpl-inter1}) and
(\ref{gpl-inter2}). As for (\ref{lem-ine3}), by Taylor expansion,
we first rewrite
\begin{eqnarray*}
r_n(\theta)&=&\frac{1}{2}\sum_{i=1}^{n}\int_{0}^1\frac
{\partial^2\log lik}
{\partial\lambda^2}(X_i;\theta,\eta_t(\theta)(Z_i))dt\{\widehat{\eta}(\theta)(Z_i)-\eta_\ast(\theta)(Z_i)\}^2
\\&\equiv&\frac{1}{2}\sum_{i=1}^{n}Q_{\theta}(X_i)\{\widehat{\eta}(\theta)(Z_i)-\eta_\ast(\theta)(Z_i)\}^2,
\end{eqnarray*}
where $\eta_t(\theta)(Z_i)=\eta_\ast(\theta)(Z_i)+t(\widehat{\eta}(\theta)
-\eta_\ast(\theta))(Z_i)$. To prove (\ref{lem-ine3}), it suffices to show that
\begin{eqnarray}
\sup_{z\in\mathcal{Z}}\left|\frac{1}{n}\sum_{i=1}^{n}\frac{\partial^j}
{\partial\theta^j}|_{\theta=\theta_0}Q_\theta(X_i)\right|
=O_P(1)\;\;\mbox{for}\;j=0,1\label{illlem}
\end{eqnarray}
in view of (\ref{pricon10}). For $j=0$, we have
$$|Q_{\theta_0}(x)|\leq\sup_{\lambda\in H}\left|\frac
{\partial^2\log lik}
{\partial\lambda^2}(x;\theta_0,\lambda)\right|=O_P(1)\;\;\mbox{for all}\;z\in\mathcal{Z}$$ based on the smoothness Condition K2. The case
$j=1$ can be established similarly. $\Box$

\begin{lemma}\label{sparlem}
Let $\eta_0(\zv)=(\eta_0(z_1),\ldots,\eta_0(z_n))'$ and $\epsv=(\epsilon_1,\ldots,\epsilon_n)'$. If $\lambda_n\rightarrow 0$, then we have
\begin{eqnarray}
\wv'A(\lambda_{n})\epsv&=&O_{P}(\lambda_{n}^{-1/(2k)}),\label{inter1a}\\
\wv'[I-A(\lambda_{n})]\eta_{0}(\zv)&=&O_{P}(n^{1/2}\lambda_n),\label{inter4a}\\
\wv'(I-A(\lambda_{n}))\wv/n&=&\Sigma+O_{P}(n^{-1/2}\vee
n^{-1}\lambda_{n}^{-1/k}).\label{inter2a}
\end{eqnarray}
\end{lemma}
{\sc Prof:}
We first state the Lemmas 4.1 and 4.3 in \citep{cw79}:
\begin{eqnarray}
&&n^{-1}\sum_{l=1}^n[(I-A(\lambda_n))\eta_0(\zv)]_{l}^{2}\leq\lambda_{n}^2J^2(\eta_0),\label{cwres0}\\
&&tr(A(\lambda_n))=O(\lambda_{n}^{-1/k}),\label{cwres}\\
&&tr(A^{2}(\lambda_n))=O(\lambda_{n}^{-1/k}).\label{cwres2}
\end{eqnarray}
%
Since
$Var[(\wv'A(\lambda_n)\epsv)_{i}]=\sigma^{2}\Sigma_{ii}tr(A^{2}(\lambda_n))$, we can show
that $[\wv'A(\lambda_n)\epsv]_{i}=O_{P}(\lambda_{n}^{-1/2k})$ based on
(\ref{cwres2}), thus proved (\ref{inter1a}). We next consider (\ref{inter4a}) by establishing that $Var[\wv'\{I-A(\lambda_n)\}\eta_0(\zv)]_i=
\Sigma_{ii}\eta_0'(\zv)[I-A(\lambda_n)]^2\eta_0(\zv)$. Then, we can prove (\ref{inter4a}) by (\ref{cwres0}).
As for (\ref{inter2a}), we first write (\ref{inter2a}) as the sum of
$$\Sigma+(\wv'\wv/n-\Sigma)-\wv'A(\lambda_n)\wv/n,$$ where the second term is $O_{P}(n^{-1/2})$ based on the central limit theorem. For the last term, we have
$E\{[\wv' A(\lambda_n)\wv]_{ij}\}^{2}=$
\begin{eqnarray*}
(\Sigma_{ij})^{2}(tr(A(\lambda_n)))^{2}+(\Sigma_{ii}\Sigma_{jj}+(\Sigma_{ij})^{2})tr(A^{2}(\lambda_n))\\+
(E(X_{1i}X_{1j})^{2}-2(\Sigma_{ij})^{2}-\Sigma_{ii}\Sigma_{jj})\sum_{r}A_{rr}^{2}(\lambda_n)
\end{eqnarray*}
for $i\neq j$. When $i=j$, we have $E|(\wv'
A(\lambda_n)\wv)_{ii}|=\Sigma_{ii}tr(A(\lambda_n))$. By considering (\ref{cwres})-(\ref{cwres2}), we have
proved (\ref{inter2a}). $\Box$

\subsection{Proof of Lemma~\ref{step}}
By
(\ref{estrel}) in Lemma~\ref{mainthm} and (\ref{res1}) in Lemma~\ref{mainthm2}, we obtain that $(\widehat{\theta}_{n}^{(k)}-\widehat{\theta}_{n})$
\begin{eqnarray*}
&=&
O_{P}\left(\frac{g_{r}\left(|t_{n}^{(k-1)}|\vee
\|\widehat{\theta}_{n}^{(k-1)}-
\widehat{\theta}_{n}\|\right)\vee nt_n^{(k-1)}\|\widehat\theta_n^{(k-1)}-\widehat\theta_n\|\vee n^{1/2-2r}}{n\{t_{n}^{(k-1)}\}^{2}}\times\right.\nonumber\\&&\mbox{\hspace{0.3in}}\|\widehat{\theta}_{n}^{(k-1)}-
\widehat{\theta}_{n}\|
\left.\vee|s_{n}^{(k-1)}|\vee\frac{g_{r}\left
(|s_{n}^{(k-1)}|\vee\|\widehat{\theta}_{n}^{(k-1)}-\widehat{\theta}_{n}\|
\right)\vee n^{1/2-2r}}
{n|s_{n}^{(k-1)}|}\right)\nonumber\\
&=&O_{P}\left(\left(|t_n^{(k-1)}|\vee\frac{n^{-2r}\vee n^{-r_{k-1}}}{|t_n^{(k-1)}|}\vee\frac{n^{-3r_{k-1}}\vee
n^{-2r-r_{k-1}}\vee n^{-1/2-2r}}{\{t_n^{(k-1)}\}^2}\right)\right.\nonumber\\
&&\mbox{\hspace{0.3in}}\left.\times n^{-r_{k-1}}\vee
\frac{n^{-3r_{k-1}}\vee n^{-2r-r_{k-1}}\vee n^{-1/2-2r}}{|s_n^{(k-1)}|}\vee|s_n^{(k-1)}|\vee
n^{-2r}\right)\nonumber\\
&=&O_{P}\left(f_{k-1}(|t_{n}^{(k-1)}|)\vee
h_{k-1}(|s_{n}^{(k-1)}|)\vee n^{-2r}\right).
\end{eqnarray*}
To analyze the above order, we have to consider three different
stages: (i) $r_{k-1}< r$; (ii) $r\leq r_{k-1}< 1/2$; (iii)
$r_{k-1}\geq1/2$. For the stage (i), the smallest order of
$f_{k-1}$, i.e., $n^{-3r_{k-1}/2}$, is achieved by taking $|t_n^{(k-1)}|\asymp n^{-r_{k-1}/2}$,
and the smallest order of $h_{k-1}$, i.e., $n^{-3r_{k-1}/2}$, is achieved by taking
$|s_n|\asymp n^{-3r_{k-1}/2}$. For the stage (ii), the
smallest order of $f_{k-1}$, i.e., $n^{-3r_{k-1}/2}$, is achieved by taking
$|t_n^{(k-1)}|\asymp n^{-r_{k-1}/2}$, and the smallest order of $h_{k-1}$, i.e.,
$n^{-(2r+r_{k-1})/2}$, is achieved by taking $|s_n^{(k-1)}|\asymp n^{-(2r+r_{k-1})/2}$.
For the last stage (iii), the
smallest order of $f_{k-1}$, i.e., $n^{-3r_{k-1}/2}$, is achieved by taking
$|t_n^{(k-1)}|\asymp n^{-r_{k-1}/2}$, and the smallest order of $h_{k-1}$, i.e.,
$n^{-r-1/4}$, is achieved by taking $|s_n^{(k-1)}|\asymp n^{-r-1/4}$. This completes the whole proof. $\Box$

\subsection{Proof of Theorem~\ref{thm1}}
According to the proof in Lemma~\ref{step}, we also need to consider the stochastic order of $\|\widehat\theta_n^{(k)}-\widehat\theta_n\|$ in terms of three stages: (i) $r_{k-1}< r$; (ii) $r\leq r_{k-1}< 1/2$; (iii) $r_{k-1}\geq1/2$. In stage (i), we have $\|\widehat{\theta}^{(k)}_n-\widehat{\theta}_n\|=O_{P}
(\|\widehat{\theta}_n^{(k-1)}-\widehat{\theta}_{n}\|^{3/2})=O_{P}(n^{-
S_1(\psi,k)})$ if $k\leq K_1(\psi,r)$. In stage (ii), we have $\|\widehat{\theta}^{(k)}_n-\widehat{\theta}_n\|=O_{P}
(\|\widehat{\theta}_n^{(k-1)}-\widehat{\theta}_{n}\|^{1/2}n^{-r})$, which implies that
$\|\widehat{\theta}_n^{(k)}-\widehat{\theta}_n\|=O_{P}(n^{-
S_2(\psi,r,k)})$ if $r\leq\psi< 1/2$. It is easy to show that $S_2(\psi,r,k)\geq1/2$ if $k\geq K_2(\psi,r,1/2)$. In the last stage (iii), we obtain the the smallest order of $\|\widehat{\theta}_n^{(k)}-\widehat{\theta}_n\|$, i.e., $O_{P}(n^{-r-1/4})$. Combining the above analysis of (i)-(iii), we can conclude that the stochastic order of $\|\widehat{\theta}^{(k)}_n-\widehat{\theta}_n\|$ is continuously improving till the optimal bound $O_P(n^{-r-1/4})$ and can be expressed as $O_P(n^{-S(\psi,r,k)})$. (\ref{kstar}) also follows from the above analysis. $\Box$

\subsection{Proof of Theorem~\ref{conv-gpl-kmle}}
We first show (\ref{gpl-mleexp}) by applying Lemma~\ref{gpl-le3}. In (\ref{gpl-kel-eq0}), we replace
$\widetilde\theta_n$ by $\widehat\theta_n$. Since $\widehat\theta_n$ is assumed to be consistent and $\theta_0$ is an interior point of $\Theta$, we have $\widehat S_n^{(1)}(\widehat\theta_n)=0$. By (\ref{lfs}) and (\ref{pricon2}), we have
\begin{eqnarray}
\;\;\;\;\;\sqrt{n}(\widehat{\theta}_{n}-\theta_0)=\sqrt{n}\widetilde{I}_0^{-1}
\mathbb{P}_n\widetilde{\ell}_0+O_P(n^{1/2-2g}\vee
n^{1/2}\|\widehat{\theta}_{n}-\theta_0\|^2)\label{gpl-int-thm}
\end{eqnarray}
given that $\widehat{\theta}_{n}$ is consistent and
$\widetilde{I}_0$ is nonsingular. Considering the range of $g$, we can show $\widehat\theta_n$ is actually $\sqrt{n}$-consistent, and thus simplify (\ref{gpl-int-thm}) to (\ref{gpl-mleexp}).

We next show (\ref{gpl-thm-1}). By (\ref{parasche-0}), we can write
$\sqrt{n}\widehat
I_n(\widehat{\theta}_n^{(0)})(\widehat{\theta}_n^{(1)}-\widehat{\theta}_n)$
as
\begin{eqnarray*}
&&\sqrt{n}\widehat
I_n(\widehat{\theta}_n^{(0)})(\widehat{\theta}_n^{(0)}-\widehat{\theta}_n)+n^{1/2}
(\widehat\ell_n(\widehat{\theta}_n^{(0)})-\widehat\ell_n(\widehat{\theta}_n))\\
&=&\sqrt{n}\widehat
I_n(\widehat{\theta}_n^{(0)})(\widehat{\theta}_n^{(0)}-
\widehat{\theta}_n)
+n^{-1/2}\widehat S_n^{(2)}(\widehat\theta_n^{(0)})(\widehat{\theta}_n^{(0)}-\widehat{\theta}_n)+
O_P(\sqrt{n}\|\widehat{\theta}_n-\widehat{\theta}_n^{(0)}
\|^2)\\
&=&O_P(\sqrt{n}\|\widehat{\theta}_n-\widehat{\theta}_n^{(0)} \|^2)
\label{intui}
\end{eqnarray*}
under Condition G. Further, by (\ref{pricon3}) and (\ref{pricon4}), we have the invertibility of
$\widehat{I}_n(\widehat{\theta}_n^{(0)})$ based on that of $\widetilde I_0$. This implies
$\widehat{\theta}_n^{(1)}-\widehat\theta_n=O_P(\|\widehat{\theta}_n^{(0)}-\widehat\theta_n\|^2)$. By the
induction principal, we can thus show
\begin{eqnarray}
\widehat{\theta}_n^{(k)}-\widehat\theta_n=
O_P(\|\widehat{\theta}_n^{(k-1)}-\widehat\theta_n\|^2)\;\;\mbox{for
any}\;k\geq 1. \label{para-inter1}
\end{eqnarray}
(\ref{gpl-thm-1}) follows from (\ref{para-inter1}) trivially.

To show (\ref{gpl-thm-2}), we first prove $\|\widehat{\theta}_{n}^{(k)}-\widehat\theta_{n}\|=$
\begin{eqnarray}
O_P\left(n^{1/2-g}\|\widehat{\theta}_{n}^{(k-1)}-\widehat{\theta}_{n}\|^2
\vee n^{-g}\|\widehat{\theta}_{n}^{(k-1)}-\widehat{\theta}_{n}\|
\right).\label{semirel2}
\end{eqnarray}
By replacing $\widetilde{\theta}_n$ and $U_n$ with $\widehat{\theta}_n$ and $(\widehat{\theta}_n^{(k-1)}-\widehat{\theta}_n)$ in
(\ref{gpl-kel-eq1}), respectively, we establish that $n^{-1/2}[\widehat S_n^{(1)}(\widehat\theta_n^{(k-1)})-\widehat S_n^{(1)}(\widehat\theta_n)]=$
\begin{eqnarray}
-\sqrt{n}\widetilde I_0(\widehat{\theta}_n^{(k-1)}
-\widehat{\theta}_n)+O_P(n^{1/2-g}\|\widehat{\theta}_n^{(k-1)}-\widehat{\theta}_n\|).\label{lem-imp}
\end{eqnarray}
Similarly, by setting $\widetilde{\theta}_n$ as
$\widehat{\theta}_n$, and then setting $U_{n}$ as
$(\widehat{\theta}_n^{(k-1)}-\widehat{\theta}_n+n^{-1/2}t_1v_j)$
and
$(\widehat{\theta}_n^{(k-1)}-\widehat{\theta}_n+n^{-1/2}t_2v_j)$ in (\ref{gpl-kel-eq1}),
respectively, we have that
\begin{eqnarray}
[\widehat{I}_n(\widehat \theta_n^{(k-1)})]_{ij}=[\widetilde I_0]_{ij}+
O_P(n^{1/2-g}\|\widehat{\theta}_n^{(k-1)}-\widehat\theta_n\| \vee
n^{-g})\label{convpi}
\end{eqnarray}
when $\widehat{I}_n^{(k-1)}$ is defined in (\ref{semicon2}). Following similar logic in analyzing (\ref{para-inter1}), we can obtain
(\ref{semirel2}) by considering (\ref{lem-imp})-(\ref{convpi}). Next we
will show that (\ref{semirel2}) implies (\ref{gpl-thm-2}) by the
following analysis. Based on (\ref{semirel2}) we have
\begin{align}\label{inter6}
\|\widehat{\theta}_{n}^{(k)}-\widehat\theta_{n}\|=
\begin{cases}
O_P(n^{-g}\|\widehat{\theta}_{n}^{(k-1)}-\widehat\theta_{n}\|) &
\text{ if  $\|\widehat{\theta}_{n}^{(k-1)}
-\widehat\theta_{n}\|=O_P(n^{-1/2})$},\\
O_P(n^{1/2-g}\|\widehat{\theta}_{n}^{(k-1)}-\widehat\theta_{n}\|^2) &
\text{ if $\|\widehat{\theta}_{n}^{(k-1)}
-\widehat\theta_{n}\|^{-1}=O_P(n^{1/2})$}.\\
\end{cases}
\end{align}
It is easy to show that
$\|\widehat{\theta}_{n}^{(L_1(\psi,g))}-\widehat{\theta}_{n}\|=O_P(n^{-1/2})$
and
$\|\widehat{\theta}_{n}^{(L_1(\psi,g)-1)}-\widehat{\theta}_{n}\|^{-1}
=O_P(n^{1/2})$. In other words, if $k\leq L_1(\psi, g)$, then we
have the relation that
$\|\widehat{\theta}_{n}^{(k)}-\widehat{\theta}_{n}\|=
O_P(n^{1/2-g}\|\widehat{\theta}_{n}^{(k-1)}
-\widehat{\theta}_{n}\|^2)$ based on (\ref{inter6}). This implies
the form of $R_1(\psi,g,k)$ in (\ref{rform}). Note that
$R_1(\psi,g,k)$ is an increasing function of $k$ under the
condition that $\psi+g>1/2$. After $L_1(\psi,g)$ iterations, we
have
\begin{eqnarray}
\|\widehat{\theta}_{n}^{(L_1(\psi,g))}-\widehat{\theta}_{n}\|
=O_P(n^{-R_1(\psi,g,L_1(\psi,g))})
=O_P(n^{-1/2}).\label{inter-para1}
\end{eqnarray}
Thus, we have the relation that
$\|\widehat{\theta}_{n}^{(k)}-\widehat{\theta}_{n}\|=O_P(n^{-g}
\|\widehat{\theta}_{n}^{(k-1)}-\widehat{\theta}_{n}\|)$ for $k\geq
(L_1(\psi,g)+1)$ based on (\ref{inter6}). Combining this relation
with (\ref{inter-para1}), we can show the form of $R_2(\psi,g,k)$
when $k>L_1(\psi,g)$. Since $R(\psi,g,k)$ is an
increasing function of $k$ given that $1/2-g<\psi\leq1/2$, the stochastic order of
$\|\widehat\theta_{n}^{(k)}-\widehat\theta_{n}\|$ is continuously
decreasing as $k\rightarrow\infty$. The calculation of $k^\ast$ also follows from the above analysis. $\Box$

\subsection{Proof of Theorem~\ref{gpl-ver}}
We first consider (\ref{pricon2}) by rewriting its LHS as
$$\frac{1}{n}\frac{\partial}
{\partial\theta}|_{\theta=\theta_0}\left[\sum_{i=1}^{n}A_{\theta,\eta_\ast(\theta)}
[\widehat{\eta}(\theta)-\eta_\ast(\theta)](X_i)+r_n(\theta)\right],$$ where $r_n(\theta)$ is defined in Lemma~\ref{gpl-le1}.
Therefore, we have
\begin{eqnarray*}
&&n^{-1}[\widehat{S}_n^{(1)}(\theta_0)-S_n^{(1)}(\theta_0)]\\
&=&\frac{1}{n}\sum_{i=1}^{n}\left(\frac{\partial}{\partial\theta}|_{\theta=\theta_0}
A_{\theta,\eta_\ast(\theta)}\right)(\widehat{\eta}_0-\eta_{\ast0})+
\frac{1}{n}\sum_{i=1}^{n}A_{\theta_0,\eta_0}(\widehat{\eta}_{0}^{(1)}
-\eta_{\ast0}^{(1)})+\frac{1}{n}\dot{r}_n(\theta_0)\\
&=&O_P(n^{-2g})
\end{eqnarray*}
by Lemma~\ref{gpl-le1} and the condition that $\delta\geq
(2g-1/2)$. As discussed previously, we will show (\ref{pricon3}) together with (\ref{pricon5}).
By Taylor expansion, we have
\begin{eqnarray*}
\widehat{S}_n(\theta)-S_n(\theta)&=&\sum_{i=1}^{n}\int_{0}^{1}
\frac{\partial\log lik}{\partial\lambda}(X_i;\theta,\eta_t(\theta)(Z_i))dt
[\widehat{\eta}(\theta)(Z_i)-\eta_\ast(\theta)(Z_i)]\\&\equiv&
\sum_{i=1}^{n}R_{\theta}(X_i)[\widehat{\eta}(\theta)(Z_i)-
\eta_\ast(\theta)(Z_i)],
\end{eqnarray*}
where $\eta_t(\theta)=\eta_\ast(\theta)+t(\widehat{\eta}(\theta)
-\eta_\ast(\theta))$. Hence, to prove (\ref{pricon3}) and (\ref{pricon5}), it suffices to
show that
\begin{eqnarray}
\sup_{\theta\in\mathcal{N}(\theta_0)}\sup_{z\in\mathcal{Z}}\left|n^{-1}\sum_{i=1}^{n}\frac{\partial^j}
{\partial\theta^j}R_{\theta}(X_i)\right|=O_P(1)\;\;\mbox{for}\;j=
0,1,2,3\label{illlem2}
\end{eqnarray}
in view of (\ref{pricon10}) and (\ref{pricon50}).
Considering the smoothness Condition K2, we can prove
(\ref{illlem2}) using the same approach as in the proof of
(\ref{illlem}).

In the end, it remains to show that the class of functions
$$\{(\partial^3/\partial\theta^3)\log
lik(x;\theta,\eta_\ast(\theta)):\theta\in\mathcal{N}(\theta_0)\}$$ is
P-Glivenko-Cantelli and that
\begin{eqnarray}
\sup_{\theta\in\mathcal{N}(\theta_0)}E\left|(\partial^3/\partial\theta^3)\log
lik(X;\theta,\eta_\ast(\theta))\right|<\infty.\label{subou}
\end{eqnarray}
Let $\ell^{(3)}(\theta,\eta(\theta))=(\partial^3/\partial\theta^3) \log
lik(x;\theta,\eta_\ast(\theta))$. For any
$\theta_1, \theta_2\in\mathcal{N}(\theta_0)$, we have $|\ell^{(3)}(\theta_1,\eta_\ast(\theta_1))-\ell^{(3)}(\theta_2,\eta_\ast(\theta_2))|$
\begin{eqnarray*}
&\leq&\sup_{\theta,\lambda}\left|\frac{\partial\ell^{(3)}}
{\partial\theta}(\theta,\lambda)\right|\|\theta_1-\theta_2\|+
\sup_{\theta,\lambda}\left|\frac{\partial\ell^{(3)}}
{\partial\lambda}(\theta,\lambda)\right|\left\|
\eta_\ast(\theta_1)-\eta_\ast(\theta_2)\right\|_\infty\\
&\leq&\sup_{\theta,\lambda}\left|\frac{\partial\ell^{(3)}}
{\partial\theta}(\theta,\lambda)\right|\|\theta_1-\theta_2\|+
\sup_{\theta,\lambda}\left|\frac{\partial\ell^{(3)}}
{\partial\lambda}(\theta,\lambda)\right|\sup_{\theta\in\mathcal{N
}(\theta_0)} \|\eta_\ast^{(1)}(\theta)\|_\infty\\&&\times\|\theta_1-\theta_2\|\\
&\leq&A\|\theta_1-\theta_2\|.
\end{eqnarray*}
By Condition K2 and
$\sup_{\theta\in\mathcal{N}(\theta_0)}\|\eta^{(1)}_\ast(\theta)\|_\infty<\infty$
in Condition C1, we know that $EA^2<\infty$. Thus, by the P-G-C
preservation Theorem 9.23 of \citep{k08} and compactness of
$\mathcal{N}(\theta_0)$, we know that
$$\{(\partial^3/\partial\theta^3)\log
lik(x;\theta,\eta_\ast(\theta)):\theta\in\mathcal{N}(\theta_0)\}$$ is
P-Glivenko-Cantelli. The last condition (\ref{subou}) follows from the Conditions K2
and C1 by some algebra. $\Box$

\subsection{Proof of Lemma~\ref{lem1}}
Let $$\widehat m_\theta(z)=\frac{\sum_{i=1}^{n}\psi_\theta(Y_i,W_i)K((z-Z_i)/b_n)}{\sum_{i=1}^n K((z-Z_i)/b_n)}.$$ Note that $\widehat\eta(\theta)(z)=\rho(\widehat m_\theta(z))$ by (\ref{inter03}). Correspondingly, we have $\eta_\ast(\theta)(z)=\rho(m_\theta(z))$ based on Lemma 7 of \cite{sw92}.
Following the proof of Lemma 8 in \citep{sw92}, we can derive that
\begin{eqnarray}
&&\sup_{\theta\in\Theta}\left\|\frac{\partial^{k+j}}{\partial z^k\partial\theta^j}
\widehat m_\theta(z)-
\frac{\partial^{k+j}}{\partial z^k\partial\theta^j}m_\theta(z)\right\|_\infty\nonumber\\&=&
O_P\left(n^{-\frac{q}{2q+4}}b_n^{-k-\frac{q+4}{q+2}}
n^{\epsilon}\vee b_n^2\right)\label{useres}
\end{eqnarray}
for any $\epsilon>0$, $k=0,1$ and $j=0,1,2,3$. Considering (\ref{inter03}), (\ref{useres}) and Condition (f), we can show that
\begin{eqnarray}
\sup_{\theta\in\mathcal{N}(\theta_0)}
\|\widehat\eta^{(s)}(\theta)-\eta^{(s)}_\ast(\theta)\|_\infty=O_P\left(n^{-\frac{q}{2q+4}}b_n^{-\frac{q+4}{q+2}}
n^{\epsilon}\vee b_n^2\right)\label{inter04}
\end{eqnarray}
for $s=0,1,2,3$ after some algebra. Following similarly logic, we show that
\begin{eqnarray}
\left\|\frac{\partial}{\partial z}\widehat\eta_0(z)-\frac{\partial}{\partial z}\eta_{\ast0}(z)\right\|_\infty&=&O_P\left(n^{-\frac{q}{2q+4}}b_n^{-\frac{2q+6}{q+2}}
n^{\epsilon}\vee b_n^2\right)\label{inter05}\\
\left\|\frac{\partial}{\partial z}\widehat\eta_0^{(1)}(z)-\frac{\partial}{\partial z}\eta_{\ast0}^{(1)}(z)\right\|_\infty&=&O_P\left(n^{-\frac{q}{2q+4}}b_n^{-\frac{2q+6}{q+2}}
n^{\epsilon}\vee b_n^2\right)\label{inter06}
\end{eqnarray}
Considering (\ref{inter04})-(\ref{inter06}), we complete the whole proof. $\Box$

\subsection{Proof of Corollary~\ref{sparseco}}
For the $\sqrt{n}$ consistency of $\widehat\theta_{\lambda_n}$, it suffices to show that, for any given $\epsilon>0$, there exists a
large constant $M$ such that
\begin{eqnarray}
P\left\{\inf_{\|s\|=M}\Delta_n(s)>0\right\}\geq 1-\epsilon,\label{inter0}
\end{eqnarray}
where $\Delta_n(s)\equiv [\widehat S_{\lambda_n}(\theta_{0}+n^{-1/2}s)-\widehat S_{\lambda_n}(\theta_{0})]$. According to (\ref{slambdan}), we have $$\Delta_n(s)\geq\widetilde S_{\lambda_n}(\theta_0+n^{-1/2}s)-\widetilde S_{\lambda_n}(\theta_0)+n\tau_n^2\sum_{j=1}^{q}\frac{|\theta_{0j}+n^{-1/2}s_j|-|\theta_{0j}|}{|\widetilde\theta_j|},$$ where $s_j$ is the $j$-th element of $s$. The Taylor expansion further gives
\begin{eqnarray}
\Delta_n(s)&\geq& n^{-1/2}s'\widetilde S_{\lambda_n}^{(1)}(\theta_0)+\frac{1}{2}s'[\widetilde S_{\lambda_n}^{(2)}(\theta_0)/n]s\nonumber\\&&+n\tau_n^2\sum_{j=1}^{q}\frac{|\theta_{0j}+n^{-1/2}s_j|-|\theta_{0j}|}{|\widetilde\theta_j|},
\label{inter7}
\end{eqnarray}
where $\widetilde S_{\lambda_n}^{(j)}(\theta_0)$ represents the $j$-th derivative of $\widetilde S_{\lambda_n}(\theta)$ at $\theta_0$. Based on (\ref{tildeslam}), we have
\begin{eqnarray}
\widetilde S_{\lambda_n}^{(1)}(\theta_0)&=&-2\wv'[I-A(\lambda_n)](\yv-\wv\theta_0),\label{inter01}\\
\widetilde S_{\lambda_n}^{(2)}(\theta_0)&=&2\wv'[I-A(\lambda_n)]\wv.\label{inter00}
\end{eqnarray}
Lemma~\ref{sparlem} implies that
\begin{eqnarray}
\widetilde S_{\lambda_n}^{(1)}(\theta_0)&=&O_P(n^{1/2}),\label{interor1}\\
\widetilde S_{\lambda_n}^{(2)}(\theta_0)&=&O_P(n)\label{interor2}
\end{eqnarray}
since $\lambda_n$ is required to converge to zero. Hence, we know the first two terms in the right hand side of (\ref{inter7}) have the same order, i.e. $O_{P}(1)$.
And the second term, which converges to some positive constant,
dominates the first one by choosing sufficiently large $M$. The
third term is bounded by $n^{1/2}\tau_n^2M_{0}$ for some
positive constant $M_{0}$ since $\widetilde{\beta}_{j}$ is the
consistent estimate for the nonzero coefficient. Considering that $\sqrt{n}\tau_{n}^2\rightarrow
0$, we have shown the $\sqrt{n}$-consistency of $\widehat\theta_{\lambda_n}$.

To complete the proof of other parts, we first need to show
\begin{eqnarray}
\|\widehat\theta_{\lambda_n}^{(1)}-\widehat\theta_{\lambda_n}\|=O_P(n^{-1})\label{asyclos}
\end{eqnarray}
based on Theorem~\ref{cor-pen}. And then we will verify Condition G' for the case $c=-2$. It is easy to show that $\mathbb{P}_n\widetilde\ell_0=\wv'\epsv/n$ and $\widetilde I_0=\Sigma$ in this example. To verify (\ref{pricon21}), we have
\begin{eqnarray*}
&&\frac{1}{n}\widehat S_{\lambda_n}^{(1)}(\theta_0)+2\mathbb{P}_n\widetilde\ell_0\\&=&-\frac{2}{n}\wv'(I-A(\lambda_n))\eta_0(\zv)+\frac{2}{n}\wv'A(\lambda_n)\epsv+\tau_n^2
\delta_n(\theta_0)\\
&=&O_P(n^{-1/2}\lambda_n\vee n^{-1}\lambda_n^{-1/(2k)}\vee\tau_n^2),
\end{eqnarray*}
where the second equality follows from Lemma~\ref{sparlem} and the fact that $\delta_n(\theta_0)=O_P(1)$. Considering the conditions on $\tau_n$ and $\lambda_n$, we have proved (\ref{pricon21}). (\ref{pricon31}) follows from (\ref{inter00}) and (\ref{inter2a}), and (\ref{pricon41}) trivially holds.
Having shown the consistency of $\widehat\theta_{\lambda_n}$ and verified G', we are able to show (\ref{asyclos}).

For any sequence of estimate $\theta_n$, the below arguments show that $\theta_n=0$ with probability tending to one if it is $\sqrt{n}$-consistent. For any $\sqrt{n}$-consistent estimator, it suffices to show that
\begin{eqnarray}
\widehat S_{\lambda_n}\{(\bar{\theta}_{1},0)\}=\min_{\|\bar{\theta}_{2}\|\leq
Cn^{-1/2}}\widehat S_{\lambda_n}\{(\bar{\theta}_{1},\bar{\theta}_{2})\}\label{inter5}
\end{eqnarray}
for any $\bar{\theta}_{1}$ satisfying
$\|\bar{\theta}_{1}-\theta_{1}\|=O_{P}(n^{-1/2})$ with
probability approaching to 1. In order to show (\ref{inter5}), we need to show
that $\partial \widehat S_{\lambda_n}(\theta)/\partial\theta_{j}<0$ for
$\theta_{j}\in(-Cn^{-1/2},0)$ and $\partial
\widehat S_{\lambda_n}(\theta)/\partial\theta_{j}>0$ for $\theta_{j}\in(0,Cn^{-1/2})$ holds when
$j=q+1,\ldots,d$ with probability tending to 1. By two term
Taylor expansion of $\widetilde S_{\lambda_n}(\theta)$ at $\theta_{0}$, $\partial
\widehat S_{\lambda_n}(\theta)/\partial\theta_{j}$ can be expressed in the following form:
\begin{eqnarray*}
\frac{\partial \widehat S_{\lambda_n}(\theta)}{\partial\theta_{j}}=\frac{\partial
\widetilde  S_{\lambda_n}(\theta_{0})}{\partial\theta_{j}}+\sum_{k=1}^{d}\frac{\partial^{2}\widetilde S_{\lambda_n}(\theta_{0})}{\partial\theta_{j}\partial\theta_{k}}(\theta_{k}-\theta_{0k})
+n\tau_{n}^2\frac{1\times
sign(\theta_{j})}{|\widetilde{\theta}_{j}|},
\end{eqnarray*}
for $j=q+1,\ldots,d$. Note that $\|\bar\theta-\theta_{0}\|=O_{P}(n^{-1/2})$ by the above
construction. Hence, we have
\begin{eqnarray*}
\frac{\partial \widehat S_{\lambda_n}(\theta)}{\partial\theta_{j}}=O_{P}(n^{1/2})+sign(\theta_{j})\frac{n\tau_{n}^2}{|\widetilde{\theta}_{j}|}
\end{eqnarray*}
by (\ref{interor1}) and (\ref{interor2}). We assume that $n^{k/(2k+1)}\tau_n\rightarrow\tau_0>0$ which implies that
$\sqrt{n}\tau_{n}^2/|\widetilde{\theta}_{j}|\rightarrow\infty$
for $\sqrt{n}$ consistent $\widetilde\theta_j$ and $j=q+1,\ldots,d$. Thus, we show that the sign of $\theta_{j}$
determines that of $\partial \widehat S_{\lambda_n}(\theta)/\partial\theta_{j}$. The above arguments apply to $\widehat\theta_{\lambda_n,2}$ and $\widehat\theta_{\lambda_n,2}^{(1)}$ since both of them are proven to be $\sqrt{n}$ consistent in view of the previous discussions, i.e., (\ref{asyclos}).

Now it remains to show the semiparametric efficiency of $\widehat\theta_{\lambda_n,1}$, which immediately implies that of $\widehat\theta_{\lambda_n,1}^{(1)}$ based on (\ref{asyclos}). Since we have shown $\widehat\theta_{\lambda_n,2}=0$, we can establish that
\begin{eqnarray}
\frac{\partial\widehat S_{\lambda_{n}}(\theta)}{\partial\theta_j}|_{\theta=(\widehat\theta_{\lambda_n,1},0)}=0\;\;\;\mbox{for any}\;j=1,\ldots,q\label{inter02}
\end{eqnarray}
with probability tending to one. Let $\wv_1$ denote the first $q$ columns of $\wv$. Applying Taylor expansion to (\ref{inter02}) around $\theta_0$, we obtain
\begin{eqnarray*}
\sqrt{n}(\widehat\theta_{\lambda_n,1}-\theta_1)&=&\sqrt{n}\left\{\frac{1}{n}\wv_1'[I-A(\lambda_n)]\wv_1\right\}^{-1}
\frac{1}{n}\wv_1'[I-A(\lambda_n)](\eta_0(\zv)+\epsv)\\
&&+O_P(\sqrt{n}\tau_n^2)\\
&=&\left\{\Sigma_{11}+O_P(n^{-1/2}\vee n^{-1}\lambda_n^{-1/k})\right\}^{-1}\frac{1}{\sqrt{n}}\wv_1'\epsv\\&&+O_P(\sqrt{n}\tau_n^2\vee n^{-1/2}\lambda_n^{-1/(2k)}\vee\lambda_n)\\
&=&\frac{1}{\sqrt{n}}\Sigma_{11}^{-1}\sum_{i=1}^{n}W_{1i}\epsilon_i+O_P(\sqrt{n}\lambda_n^{2}\vee\sqrt{n}\tau_n^2)
\end{eqnarray*}
based on (\ref{inter01}) \& (\ref{inter00}). This completes the whole proof. $\Box$

\subsection{Proof of Theorem~\ref{init}} Define
$\mathcal{N}_{n}=\{\theta:\|\theta-\theta_{0}\|\leq Mn^{-\psi}\}$
and $\mathcal{N}_{n}^{c}$ as its complement for any $0<M<\infty$. Note that
$\mathcal{D}_{n}\cap\mathcal{N}_{n}\neq\emptyset$ for large enough $M$ and $\mathcal{D}_{n}\cap\mathcal{N}_{n}^{c}\neq\emptyset$
for large enough $n$. We first consider (\ref{thm5}). For sufficiently large $M$ and any $C_1>0$, we have
\begin{eqnarray*}
P\left(\theta_{n}^{D}\in\mathcal{N}_{n}^{c}\right)&=&P\left(\theta_n^D\in \mathcal{N}_n^c\;\;\mbox{and}\;\;\theta_{iD}\in\mathcal{N}_n\;\;\mbox{for some}\;i\right)\\
&\leq&
P\left(\max_{\mathcal{D}_{n}\cap\mathcal{N}_{n}}\widehat S_n(\theta)
\leq\max_{\mathcal{D}_{n}\cap\mathcal{N}_{n}^{c}}\widehat S_n(\theta)\right)\\
&\leq&
P\left(\max_{\mathcal{D}_{n}\cap\mathcal{N}_{n}}\widehat S_n(\theta)<
\widehat S_{n}(\theta_{0})-C_{1}n^{1-2\psi}\right)\\
&&+
P\left(\left\{\max_{\mathcal{D}_{n}\cap\mathcal{N}_{n}}\widehat S_n(\theta)
\leq\max_{\mathcal{D}_{n}\cap\mathcal{N}_{n}^{c}}\widehat S_n(\theta)\right\}\right.\\
&&\;\;\;\;\;\;\;\;\;\;\;\;\;\;\;\;\;\;\;\;\;\;\;\;\;\;\;\;\cap\left.\left\{\max_{\mathcal{D}_{n}\cap\mathcal{N}_{n}}\widehat S_n(\theta)\geq
\widehat S_n(\theta_{0})-C_{1}n^{1-2\psi}\right\}\right)\\
&\leq&P\left(\max_{\mathcal{D}_{n}\cap\mathcal{N}_{n}}n^{-1/2}(\widehat S_n(\theta)-\widehat S_n(\theta_{0}))<-C_{1}n^{1/2-2\psi}\right.\\
&&\;\;\;\;\;\;\;\;\;\;\;\;\;\;\;\;\;\;\;\;\;\;\;\;\;\;\;\;\;\;\;\;\;
\;\;\;\;\;\;\;\;\;\;\;\;\;\;\;\;\;\;\;\;\;\;\;\;\;\cap\left.\{\theta_n^o\;\mbox{is consistent}\}\right)\\
&&+P\left(\max_{\mathcal{N}_{n}^{c}}n^{-1/2}(\widehat S_n(\theta)-\widehat S_n(\theta_{0}))
\geq-C_{1}n^{1/2-2\psi}\right)\\&&+P\left(\theta_n^o\;\mbox{is inconsistent}\right)\\
&\leq&I+II+III,
\end{eqnarray*}
where $\theta_n^o=\arg\max_{\mathcal{D}_n\cap\mathcal{N}_n}\widehat S_n(\theta)$.

The definition of $\mathcal{N}_n$ implies $III\rightarrow 0$ for any $M$ as $n\rightarrow\infty$.
We next analyze the term I as follows. In view of (\ref{iniass1}) and the definition of $\mathcal{N}_n$, we have that
\begin{eqnarray*}
I&=&P\left(\sqrt{n}(\theta_n^o-\theta_{0})'\mathbb{P}_{n}\widetilde{\ell}_{0}
-\frac{\sqrt{n}}{2}(\theta_n^o-\theta_{0})'\widetilde{I}_{0}(\theta_n^o-\theta_{0})+n^{-1/2}\Delta_{n}(\theta_n^o)
\right.\\
&&\;\;\;\;\;\;\;\;\;\;\;\;\;\;\;\;\;\;\;\;\;\;\;\;\;\;\;\;\;\;\;\;\;
\;\;\;\;\;\;\;\;\;\;\;\;\;\;\;\;\;\;\;\;\;\;\;\;\;\;\;\;\;\;\;\;\;\;\;\;\;\;\;\;\;\;\;\;\;\;\;\left.<-C_{1}n^{1/2-2\psi}\right)\nonumber\\
&\leq&
P\left(\|\sqrt{n}\mathbb{P}_{n}\widetilde{\ell}_{0}\|
\|\theta_n^o-\theta_{0}\|+(\delta_{max}\sqrt{n}/2)\|\theta_n^o-\theta_{0}\|^2+\|n^{-1/2}\Delta_{n}(\theta_n^o)\|\right.\\
&&\;\;\;\;\;\;\;\;\;\;\;\;\;\;\;\;\;\;\;\;\;\;\;\;\;\;\;\;\;\;\;\;\;
\;\;\;\;\;\;\;\;\;\;\;\;\;\;\;\;\;\;\;\;\;\;\;\;\;\;\;\;\;\;\;\;\;\;\;\;\;\;\;\;\;\;\;\;\;\;\;\;\;\;\left.>C_{1}n^{1/2-2\psi}\right)\nonumber\\
&\leq& P\left(\|\sqrt{n}\mathbb{P}_{n}\widetilde{\ell}_{0}\|>
\frac{C_{1}-\delta_{max} M^{2}/2}{M}n^{1/2-\psi}+o_{P}(n^{1/2-\psi})\right)\\
&\leq&\bar{I},
\end{eqnarray*}
where $\delta_{max}$ is the largest eigenvalue of $\widetilde{I}_{0}$, and the second inequality follows from the definitions of $\mathcal{N}_n$ and $\Delta_{n}$, and the range that $2r>1/2\geq\psi>0$. Denote $\theta_n^\ast=\arg\max_{\mathcal{N}_n^c}\widehat S_n(\theta)$.
We will show $II\rightarrow 0$ by first decomposing it as $II_1+II_2$, where
\begin{eqnarray*}
II_1&=&P\left(n^{-1/2}(\widehat S_n(\theta_n^\ast)-\widehat S_n(\theta_{0}))
\geq-C_{1}n^{1/2-2\psi}\cap\{\theta_n^\ast\;\mbox{is consistent}\}\right),\\
II_2&=&P\left((\widehat S_n(\theta_n^\ast)-\widehat S_n(\theta_{0}))
\geq-C_{1}n^{1-2\psi}\cap\{\theta_n^\ast\;\mbox{is inconsistent}\}\right).
\end{eqnarray*}
Note that we can write  $n^{-1/2}\Delta_{n}(\theta_n^\ast)$ as $\sqrt{n}\|\theta_n^\ast-\theta_0\|^2\epsilon_{1n}+\sqrt{n}\|\theta_n^\ast-\theta_0\|\epsilon_{2n}$, where $\epsilon_{1n}=o_P(1)$ and $\epsilon_{2n}=o_P(n^{-1/2})$, in the event that $\{\theta_n^\ast\;\mbox{is consistent}\}$. Thus, according to (\ref{iniass1}), we can write $II_1$ as
\begin{eqnarray*}
&&P\left((\theta_{n}^{\ast}-\theta_{0})'\sqrt{n}\mathbb{P}_{n}\widetilde{\ell}_{0}+\sqrt{n}\|\theta_n^\ast-\theta_0\|\epsilon_{2n}\geq\frac{\sqrt{n}}{2}
(\theta_{n}^{\ast}-\theta_{0})'\widetilde{I}_{0}(\theta_{n}^{\ast}-\theta_{0})\right.\nonumber\\
&&\;\;\;\;\;\;\;\;\;\;\;\;\;\;\;\;\;\;\;\;\;
\;\;\;\;\;\;\;\;\;\;\;\;\;\;\;\;\;\;\;\;\;\;\;\;\;\;\;\;\;\;\;\left.-\sqrt{n}\|\theta_n^\ast-\theta_0\|^2\epsilon_{1n}
-C_{1}n^{1/2-2\psi}\right)\\
&\leq&P\left(\|\theta_{n}^{\ast}-\theta_{0}\|\left[\|\sqrt{n}\mathbb{P}_{n}\widetilde{\ell}_{0}\|+\sqrt{n}\epsilon_{2n}\right]\geq
\frac{\sqrt{n}}{2}\|\theta_n^\ast-\theta_0\|^2\delta_{min}\right.\nonumber\\
&&\;\;\;\;\;\;\;\;\;\;\;\;\;\;\;\;\;\;\;\;\;
\;\;\;\;\;\;\;\;\;\;\;\;\;\;\;\;\;\;\;\;\;\;\;\;\;\;\;\;\;\;\;\left.-\sqrt{n}\|\theta_n^\ast-\theta_0\|^2\epsilon_{1n}
-C_{1}n^{1/2-2\psi}\right)\\
&\leq&P\left(\left[\|\sqrt{n}\mathbb{P}_{n}\widetilde{\ell}_{0}\|+\sqrt{n}\epsilon_{2n}\right]\geq\sqrt{n}\|\theta_n^\ast-\theta_0\|(\delta_{min}/2-
\epsilon_{1n})-\frac{C_1n^{1/2-\psi}}{K}\right)\\
&\leq&P\left(\left[\|\sqrt{n}\mathbb{P}_{n}\widetilde{\ell}_{0}\|+\sqrt{n}\epsilon_{2n}\right]\geq\frac{\delta_{min}K^2/2-C_1}{K}n^{1/2-\psi}+o_P(n^{1/2-\psi})\right)\\
&\leq& \bar{II}_{1},
\end{eqnarray*}
where $\delta_{min}>0$ is the smalest eigenvalue of $\widetilde I_0$. All the above inequalities follow from the fact that $\|\theta_{n}^\ast-\theta_{0}\|\geq Kn^{-\psi}$ for some $K>M$ and $\epsilon_{1n}=o_P(1)$. The term $II_2$ is shown to converge to zero by the following contradiction arguments. By assuming that the event $\{(\widehat S_n(\theta_n^D)-\widehat S_n(\theta_{0}))
\geq-C_{1}n^{1-2\psi}\}$ holds, we have $|\widehat S_n(\theta_n^D)-\widehat S_n(\widehat\theta_n)|=\widehat S_n(\widehat\theta_n)-\widehat S_n(\theta_n^D)\leq \widehat S_n(\widehat\theta_n)-\widehat S_n(\theta_0)+C_1n^{1-2\psi}$. Note that (\ref{iniass1}) and the consistency of $\widehat\theta_n$ implies $\widehat S_n(\theta_0)-\widehat S_n(\widehat\theta_n)=o_P(n)$. Then, we can show that $|\widehat S_n(\theta_n^D)-\widehat S_n(\widehat\theta_n)|/n=o_P(1)$ which implies that $\theta_n^D$ is consistent by (\ref{asyuni}). This implication contradicts with another event in $II_2$, i.e., $\{\theta_n^D\;\mbox{is inconsistent}\}$. Therefore we can claim that $II_2\rightarrow 0$.

In view of the above discussions, it remains to show that $\bar I$ and $\bar {II}_1$ converge to zero.
Note that $\|\sqrt{n}\mathbb{P}_{n}\widetilde{\ell}_{0}\|$ in $\bar I$ is $O_P(1)$, and so is $(\|\sqrt{n}\mathbb{P}_{n}\widetilde{\ell}_{0}\|+\sqrt{n}\epsilon_{2n})$ in $\bar{II}_1$. Therefore, by choosing sufficiently large $C_1$ and $K>M$, meanwhile keeping the inequality $\delta_{max}M^2<2C_1<\delta_{min}K^2$ valid, we show that $\bar I$ and $\bar {II}_1$ can be arbitrarily close to zero. For example, we can take $K=M+B$ and $C_1=(\delta_{max}M^2+\delta_{min}(M+B)^2)/4$ for some fixed $B>0$ and sufficiently large $M$. This completes the proof of (\ref{thm5}).

Our proof of (\ref{thm6}) is similar as that of (\ref{thm5}). Denote $\theta_{iS}$ as an element in $\mathcal{S}_n$. Similarly, we have
\begin{eqnarray*}
P(\theta_{n}^{S}\in\mathcal{N}_{n}^{c})
&\leq&E\left\{P\left(\theta_n^S\in \mathcal{N}_n^c\;\;\mbox{and}\;\;\theta_{iS}\in\mathcal{N}_n\;\;\mbox{for some}\;i|\mathcal{S}_n\right)\right\}\\&&+E\left\{P\left(
\theta_{iS}\in\mathcal{N}_n^c\;\;\mbox{for all}\;i|\mathcal{S}_n\right)\right\}\\
&\leq&P\left(\max_{\mathcal{S}_{n}\cap\mathcal{N}_{n}}\widehat S_n(\theta)
\leq\max_{\mathcal{S}_{n}\cap\mathcal{N}_{n}^{c}}\widehat S_n(\theta)\right)+P\left(
\theta_{iS}\in\mathcal{N}_n^c\;\;\mbox{for all}\;i\right)\\
&\leq&P\left(\max_{\mathcal{S}_{n}\cap\mathcal{N}_{n}}n^{-1/2}(\widehat S_n(\theta)-\widehat S_n(\theta_{0}))<-C_{2}n^{1/2-2\psi}\right)\\
&&+P\left(\max_{\mathcal{S}_n\cap\mathcal{N}_{n}^{c}}n^{-1/2}(\widehat S_n(\theta)-\widehat S_n(\theta_{0}))\geq-C_{2}n^{1/2-2\psi}
\right)\\&&+P\left(
\theta_{iS}\in\mathcal{N}_n^c\;\;\mbox{for all}\;i\right)\\
&\leq&P\left(\max_{\mathcal{S}_{n}\cap\mathcal{N}_{n}}n^{-1/2}(\widehat S_n(\theta)-\widehat S_n(\theta_{0}))<-C_{2}n^{1/2-2\psi}\right.\nonumber\\
&&\;\;\;\;\;\;\;\;\;\;\;\;\;\;\;\;\;\;\;\;\;
\;\;\;\;\;\;\;\;\;\;\;\;\;\;\;\;\;\;\;\;\;\;\;\;\;\;\;\;\;\;\;\;\;\;\;\;\left.\cap
\{\theta_n^\dag\;\mbox{is consistent}\}\right) \\&&
+P\left(\max_{\mathcal{N}_{n}^{c}}n^{-1/2}(\widehat S_n(\theta)-\widehat S_n(\theta_{0}))\geq-C_{2}n^{1/2-2\psi}\right)\\&&
+P(\theta_n^\dag\;\mbox{is inconsistent})+P\left(
\theta_{iS}\in\mathcal{N}_n^c\;\;\mbox{for all}\;i\right)\\
&\leq&I'+II'+III'+IV',
\end{eqnarray*}
where $C_2$ is an arbitrary positive constant and $\theta_n^\dag=\arg\max_{\mathcal{S}_n\cap\mathcal{N}_n}\widehat S_n(\theta)$.

We first consider the terms $III'$ \& $IV'$. Since $\theta_n^\dag\in\mathcal{N}_n$, we have $III'\rightarrow 0$ for any $M$ as $n\rightarrow\infty$. The term $IV'$ is computed as
\begin{eqnarray}
\left(1-P(\bar\theta\in\mathcal{N}_n)\right)^{card(\mathcal{S}_n)}.\label{ivterm}
\end{eqnarray}
Since the density of $\bar\theta$ is assumed to be bounded away from zero around $\theta_0$ and $card(\mathcal{S}_n)\geq\widetilde Cn^{\psi}$, (\ref{ivterm}) is bounded above by
\begin{eqnarray}
\left(1-\rho n^{-\psi}M\right)^{card(\mathcal{S}_n)}&\leq&\left(1-\rho M\widetilde C/card(\mathcal{S}_n)\right)^{card(\mathcal{S}_n)}\nonumber\\&\longrightarrow&\exp(-\rho M\widetilde C),\label{ivterm2}
\end{eqnarray}
for some $\rho>0$.

We next consider $I'$. According to (\ref{iniass1}), we can show
\begin{eqnarray*}
&&n^{-1/2}(\widehat S_n(\theta_n^\dag)-\widehat S_n(\theta_0))\\&\geq&\max_{\mathcal{S}_n\cap\mathcal{N}_n}\left\{-\frac{\sqrt{n}}{2}(\theta-\theta_0)'
\widetilde I_0(\theta-\theta_0)\right\}-\max_{\mathcal{S}_n\cap\mathcal{N}_n}\{-\sqrt{n}(\theta-\theta_0)'\mathbb{P}_n\widetilde\ell_0
-\Delta_{n}(\theta)/\sqrt{n}\}\\
&\geq&-\min_{\mathcal{S}_n\cap\mathcal{N}_n}\left\{\frac{\sqrt{n}}{2}(\theta-\theta_0)'
\widetilde I_0(\theta-\theta_0)\right\}-\max_{\mathcal{S}_n\cap\mathcal{N}_n}\{-\sqrt{n}(\theta-\theta_0)'\mathbb{P}_n\widetilde\ell_0
-\Delta_{n}(\theta)/\sqrt{n}\}.
\end{eqnarray*}
Therefore, we can bound $I'$ by $I_1'+I_2'$, where
\begin{eqnarray*}
I_1'&=&P\left(\max_{\mathcal{S}_n\cap\mathcal{N}_n}\{-\sqrt{n}(\theta-\theta_0)'\mathbb{P}_n\widetilde\ell_0-n^{-1/2}\Delta_{n}(\theta)\}>
(C_2/2)n^{1/2-2\psi}\right),\\
I_2'&=&P\left(\min_{\mathcal{S}_n\cap\mathcal{N}_n}\{\sqrt{n}(\theta-\theta_0)'\widetilde I_0(\theta-\theta_0)\}>C_2n^{1/2-2\psi}\right).
\end{eqnarray*}
Given sufficiently large $C_2/M$, $I_1'$ can be arbitrarily close to zero since
\begin{eqnarray}
I_1'&\leq&P\left(\|\sqrt{n}\mathbb{P}_n\widetilde\ell_0\|>\frac{C_2}{2M}n^{1/2-\psi}+O_P(n^{1/2-2\psi}\vee n^{1/2-2r})\right)
\nonumber\\
&\leq&P\left(\|\sqrt{n}\mathbb{P}_n\widetilde\ell_0\|>\frac{C_2}{2M}n^{1/2-\psi}+o_P(n^{1/2-\psi})\right),\label{term1}
\end{eqnarray}
where the last inequality follows from the assumption that $2r>1/2\geq\psi$. Since $\min_{\mathcal{N}_n^c}\{\sqrt{n}(\theta-\theta_0)'\widetilde I_0(\theta-\theta_0)\}>C_2n^{1/2-2\psi}$ by choosing $\delta_{\min}M^2>C_2$, $I_2'$ is bounded above by
\begin{eqnarray}
&&P\left(\min_{\mathcal{S}_n}\{\sqrt{n}(\theta-\theta_0)'\widetilde I_0(\theta-\theta_0)\}>C_2n^{1/2-2\psi}\right)\nonumber\\
&\leq&\left[P(\sqrt{n}(\bar\theta-\theta_0)'\widetilde I_0(\bar\theta-\theta_0)>C_2n^{1/2-2\psi})\right]^{card(\mathcal{S}_n)}\nonumber\\
&\leq&\left[1-P(\|\bar\theta-\theta_0\|\leq(C_2/\delta_{max})^{1/2}n^{-\psi})\right]^{card(\mathcal{S}_n)}\nonumber\\
&\leq&\left[1-P\left(\|\bar{\theta}-\theta_{0}\|\leq (C_2/\delta_{max})^{1/2}\widetilde C/card(\mathcal{S}_n)\right)\right]^{card(\mathcal{S}_{n})}\nonumber\\
&\leq&(1-\rho\widetilde C(C_{2}/\delta_{max})^{1/2}/card(\mathcal{S}_{n}))^{card(\mathcal{S}_{n})}\nonumber\\&\longrightarrow&\exp(-\rho\widetilde C \sqrt{C_{2}/\delta_{max}})\label{term12}
\end{eqnarray}
for some $\rho>0$. In the above, the third and fourth inequality follows from the assumptions that $card(\mathcal{S}_n)\geq\widetilde Cn^{\psi}$ and the density for $\bar\theta$ is bounded away from zero around $\theta_0$, respectively. By assuming that $2C_2<K^2\delta_{min}$ for some $K>M$,
we can prove that $II'\rightarrow 0$ in the same manner as we show $II\rightarrow 0$.
%
%

Let $L=\min\{K^2/2, M^2\}$. In view of (\ref{ivterm2}), (\ref{term1}), (\ref{term12}) and the above discussions on $II'$, by choosing sufficiently large $C_2$, $K>M$ and $C_2/M$, meanwhile keeping the inequality $C_2<L\delta_{min}$ valid, we can make $P(\theta_n^{S}\in\mathcal{N}_n^c)$ arbitrarily small. For example, we can take $C_2=M^{3/2}\delta_{min}$ and $K=M+B$, for some fixed $B>0$ and sufficiently large $M$. This completes the whole proof. $\Box$

\end{document}